\documentclass[a4paper]{article}
\usepackage{amssymb}
\usepackage{mathrsfs}
\usepackage{amsmath}
\usepackage{amsfonts,amssymb}
\catcode`\@=11 \@addtoreset{equation}{section}

\catcode`\@=12

\newtheorem{Theorem}{Theorem}[section]
\newtheorem{Definition}{Definition}[section]
\newtheorem{Lemma}{Lemma}[section]

\newtheorem{Remark}{Remark}[section]

\newcommand{\p}{\partial}

\def\2{{I \hskip -1.0mm I}}
\def\3{{I \hskip -1.0mm I\hskip -1.0mm I}}
\def\4{{I \hskip -0.9mm V}}
\def\6{{V \hskip -1.35mm I}}

\title{Hyperbolic geometric flow (I): short-time existence and
nonlinear stability}
\author{Wen-Rong Dai$^1$, \;\; De-Xing Kong$^{1,}\footnote{Corresponding author: kong@cms.zju.edu.cn.}$
\;\; and \;\; Kefeng Liu$^{2,1}$\\
\\ $^{1}$Center of Mathematical
Sciences, Zhejiang University\\
 Hangzhou 310027, China\\
$^2$Department of Mathematics, UCLA, CA 90095, USA}
\date{ }
\begin{document}
\maketitle
\begin{abstract} In this paper we establish the short-time existence
and uniqueness theorem for hyperbolic geometric flow, and prove the
nonlinear stability of hyperbolic geometric flow defined on the
Euclidean space with dimension larger than 4. Wave equations
satisfied by the curvatures are derived. The relations of the
hyperbolic geometric flow with the Einstein equations and the Ricci
flow are discussed.

\vskip 6mm


\vskip 3mm

\end{abstract}
\newpage
\baselineskip=7mm

\section{Introduction}

This is the first of a series of papers devoted to the study of
hyperbolic geometric flow and its applications to geometry and
physics. Hyperbolic geometric flow was first studied by Kong and
Liu in \cite{kl}. To introduce such flow we were partially
motivated by the Einstein equations in general relativity and the
recent progress in the Hamilton's Ricci flow, and by the
possibility of applying the powerful theory of hyperbolic partial
differential equations to geometry. Hyperbolic geometric flow is a
system of nonlinear evolution partial differential equations of
second order, it is very natural to understand certain wave
phenomena in nature as well as the geometry of manifolds, in
particular, it describes the wave character of the metrics and
curvatures of manifolds. We will see that the hyperbolic geometric
flow carries many interesting features of both the Ricci flow as
well as the Einstein equations. It has many promising applications
to both subjects.

The elliptic and parabolic partial differential equations have been
successfully applied to differential geometry and physics. Typical
examples are the Hamilton's Ricci flow and Schoen-Yau's solution of
the positive mass conjecture. A natural and important question is if
we can apply the well-developed theory of hyperbolic differential
equations to solve problems in differential geometry and theoretical
physics. This series of papers is an attempt to apply the hyperbolic
equation techniques to study some geometrical problems and physical
problems. One has found interesting results in these directions, see
for example
%\cite{ckl1}-\cite{ckl2} and
\cite{shu} for the applications of the hyperbolic geometric flow
equations to physics. Our results already show that the hyperbolic
geometric flow is a natural and powerful tool to study some
important problems arising from differential geometry and general
relativity such as singularities, existence and regularity. In this
paper we study the basic properties of the hyperbolic geometric flow
such as the short-time existence, nonlinear stability and the wave
feature of the curvatures. In the sequel we will study several
fundamental problems, for example, long-time existence, formation of
singularities, as well as the physical and geometrical applications.

Let $\mathscr{M}$  be an $n$-dimensional complete Riemannian
manifold with Riemannian metric $g_{ij}$, the Levi-Civita connection
is given by the Christoffel symbols
\begin{equation*}
\Gamma^{k}_{ij}=\frac{1}{2}g^{kl}\left\{\frac{\partial
g_{jl}}{\partial x^{i}}+\frac{\partial g_{il}}{\partial
x^{j}}-\frac{\partial g_{ij}}{\partial x^{l}}\right\},
\end{equation*}
where $g^{ij}$ is the inverse of $g_{ij}$. The Riemannian curvature
tensors read
\begin{eqnarray*}
R^{k}_{ijl}=\frac{\partial\Gamma^{k}_{jl}}{\partial
x^{i}}-\frac{\partial\Gamma^{k}_{il}}{\partial
x^{j}}+\Gamma^{k}_{ip}\Gamma^{p}_{jl}-\Gamma^{k}_{jp}\Gamma^{p}_{il},\;\;\;
R_{ijkl}=g_{kp}R^{p}_{ijl}.
\end{eqnarray*}
The Ricci tensor is the contraction
\begin{equation*}
R_{ik}=g^{jl}R_{ijkl}
\end{equation*}
and the scalar curvature is
\begin{equation*}
R=g^{ij}R_{ij}.
\end{equation*}
The hyperbolic geometric flow under the consideration is the
following evolution equation
\begin{equation}
\frac{\partial^{2}g_{ij}}{\partial t^{2}}=-2R_{ij}
\end{equation}
for a family of Riemannian metrics $g_{ij}(t)$ on $\mathscr{M}$.
More general hyperbolic geometric flows were also introduced in
\cite{kl}. A natural and fundamental problem is the short-time
existence and uniqueness theorem of hyperbolic geometric flow (1.1).
In the present paper, we prove the following short-time existence
and uniqueness theorem, the nonlinear stability theorem for
Euclidean space, and derive the corresponding wave equations for the
curvatures. These results were announced in Kong and Liu \cite{kl}.
\begin{Theorem}
Let $(\mathscr{M},g^0_{ij}(x))$ be a compact Riemannian manifold.
Then there exists a constant $h>0$ such that the initial value
problem
$$
\begin{cases}
\frac{\partial^2}{\partial t^2}g_{ij}(x,t)=-2R_{ij}(x,t),\\
g_{ij}(x,0)=g^0_{ij}(x),\ \frac{\partial g_{ij}}{\partial
t}(x,0)=k^0_{ij}(x),
\end{cases}$$
has a unique smooth solution $g_{ij}(x,t)$ on $\mathscr{M}\times
[0,h]$, where $k^0_{ij}(x)$ is a symmetric tensor on
$\mathscr{M}$.\end{Theorem}

The main difficulty to prove this theorem is that, the hyperbolic
geometric flow (1.1) is a system of nonlinear weakly-hyperbolic
partial differential equations of second order. The degeneracy of
the system is caused by the diffeomorphism group of $\mathscr{M}$
which acts as the gauge group of the hyperbolic geometric flow.
Because the hyperbolic geometric flow (1.1) is only weakly
hyperbolic, the short-time existence and uniqueness result on a
compact manifold does not come from the standard PDEs theory
directly. In order to prove the above short-time existence and
uniqueness theorem, using the gauge fixing idea as in the Ricci
flow, we can derive a system of nonlinear strictly-hyperbolic
partial differential equations of second order, thus Theorem 1.1
comes from the standard PDEs theory. On the other hand, we can
reduce the hyperbolic geometric flow (1.1) to a quasilinear
symmetric hyperbolic system of first order, then using the
Friedrich's theory \cite{f} of symmetric hyperbolic system (more
exactly, the quasilinear version \cite{fm}) we can also prove
Theorem 1.1.

Noting an important result on nonlinear wave equations (see
\cite{k}), we will find its other interesting application to
geometry, by applying this result to the wave equations of
curvatures. In the present paper we first use it to prove the
nonlinear stability of the flat solution of the hyperbolic geometric
flow defined on the Euclidean space with dimension larger than 4.
More precisely, we have
\begin{Theorem}
The flat metric $g_{ij}=\delta_{ij}$ of the Euclidean space
$\mathbb{R}^{n}$ with $n\geq 5$ is nonlinearly stable.
\end{Theorem}

See Section 4 for the precise definition of nonlinear stability
which is very important in general relativity. The key point of the
proof of this theorem is the global existence of classical solutions
of the Cauchy problem for the nonlinear wave equations.

Similar to Hamilton \cite{h}, we derive the corresponding wave
equations for the curvatures, for example, we have
\begin{Theorem} Under the hyperbolic geometric flow (1.1), the
curvature tensor satisfies the evolution equation
\begin{eqnarray}\begin{array}{lll}
\dfrac{\partial^2}{\partial t^2}R_{ijkl}= & \triangle
R_{ijkl}+2\left(B_{ijkl}-B_{ijlk}-B_{iljk}+B_{ikjl}\right)\vspace{2mm}\\
&
-g^{pq}\left(R_{pjkl}R_{qi}+R_{ipkl}R_{qj}+R_{ijpl}R_{qk}+R_{ijkp}R_{ql}\right)\vspace{2mm}\\
& +2g_{pq}\left(\dfrac{\partial\Gamma^{p}_{il}}{\partial t}
\dfrac{\partial\Gamma^{q}_{jk}}{\partial
t}-\dfrac{\partial\Gamma^{p}_{jl}}{\partial t}
\dfrac{\partial\Gamma^{q}_{ik}}{\partial t}\right),\end{array}
\end{eqnarray}
where $B_{ijkl}=g^{pr}g^{qs}R_{piqj}R_{rksl}$ and $\triangle$ is the
Laplacian with respect to the evolving metric.
\end{Theorem}

The wave equations for the Ricci and scalar curvatures are stated
and proved in Section 5. This is similar to the Ricci flow equation,
except that there are quadratic lower order terms involving the
connection coefficients. It turns out that there is a very rich
theory in nonlinear wave equations to deal with such terms, see
\cite{k}.

From the above results we can already see that the hyperbolic
geometric flow has many features of the Ricci flow, therefore many
well-developed techniques in the Ricci flow may be applied to the
study of this new kind of flow equations. On the other hand, the
hyperbolic geometric flow can also be viewed as the leading terms of
the vacuum Einstein equations. Since the hyperbolic geometric flow
contains the major terms in the Einstein equations, it not only
becomes simpler and more symmetric, but also possesses rich and
beautiful geometric properties. See Section 6 for more detailed
discussions on the relations between the hyperbolic geometric flow
and the Einstein equations, and more generally its relations with
other important problems in general relativity.

The paper is organized as follows. In Section 2, using the gauge
fixing idea as in the Ricci flow, we derive a system of nonlinear
strictly-hyperbolic partial differential equations of second order.
In Section 3 we reduce the hyperbolic geometric flow (1.1) to a
quasilinear symmetric hyperbolic system of first order, and give the
proof of Theorem 1.1. Section 4 is devoted to the nonlinear
stability of the hyperbolic geometric flow defined on the Euclidean
space with dimension larger than 4. In Section 5, we derive the wave
equations satisfied by the curvatures, and illustrate the wave
character of the curvatures. Some discussions are given in Section
6.

\section{Strict hyperbolicity of hyperbolic geometric flow}

In this section we consider a modified system of evolution
equations of the hyperbolic geometric flow, which is strictly
hyperbolic so that we can get a solution for a short time by
solving the corresponding Cauchy problem. The solution of the
system (1.1) then comes from the solution of the modified
equations.

Let $\mathscr{M}$ be a compact $n$-dimensional manifold. We
consider the hyperbolic geometric flow (1.1) on $\mathscr{M}$,
that is,
\begin{equation}
\frac{\partial^{2}}{\partial t^{2}}{g_{ij}}(x,t) =
-{2}{R}_{ij}(x,t).
\end{equation}
Suppose ${\hat{g}_{ij}(x,t)}$  is a solution of the hyperbolic
geometric flow (2.1), and ${\psi_{t}}$ : $\mathscr{M}\rightarrow
\mathscr{M}$ is a family of diffeomorphisms of $\mathscr{M}$. Let
$$g_{ij}(x,t) = {\psi_{t}^\ast\hat{g}_{ij}}(x,t)$$
be the pull-back metrics. We now want to find the evolution
equations for the metrics $g_{ij}(x,t)$. Denote by $y(x,t) =
{\varphi_{t}(x)} =(y^1(x,t),y^{2}(x,t),\cdots,y^{n}(x,t))$ in
local coordinates. Then
\begin{equation}
g_{ij}(x,t) = \frac{\partial y^\alpha}{\partial x^{i}}
\frac{\partial y^\beta}{\partial x^{j}} \hat{g}_{\alpha\beta}(y,t)
\end{equation}
and
\begin{eqnarray}\begin{array}{lll}
\dfrac{\partial}{\partial t}{g_{ij}(x,t)} &=&
\dfrac{\partial}{\partial t}
\left[{\hat{g}_{\alpha\beta}(y,t)}\dfrac{\partial y^\alpha}{\partial
x^{i}} \dfrac{\partial y^\beta}{\partial x^{j}}\right]\vspace{2mm}
\\&=& \dfrac{\partial y^\alpha}{\partial x^{i}}\dfrac{\partial
y^\beta}{\partial x^{j}}\dfrac{d}{d
t}\hat{g}_{\alpha\beta}(y(x,t),t)+\hat{g}_{\alpha\beta}(y,t)\dfrac{\partial}{\partial
t}\left(\dfrac{\partial y^\alpha}{\partial x^{i}}\dfrac{\partial
y^\beta}{\partial x^{j}}\right).\nonumber\end{array}\end{eqnarray}
Furthermore, we have
\begin{eqnarray}\begin{array}{lll}
\dfrac{\partial^{2}}{\partial t^{2}}{g_{ij}}(x,t) &=&
\dfrac{\partial y^\alpha}{\partial x^{i}} \dfrac{\partial
y^\beta}{\partial x^{j}} \dfrac{d^{2}\hat{g}_{\alpha\beta}}{d
t^{2}}(y(x,t),t)+\dfrac{\partial}{\partial
x^{i}}\left(\dfrac{\partial^{2} y^\alpha}{\partial
t^{2}}\right)\dfrac{\partial y^\beta}{\partial
x^{j}}{\hat{g}_{\alpha\beta}}\vspace{2mm}\\& & +\dfrac{\partial
y^\alpha}{\partial x^{i}}\dfrac{\partial}{\partial
x^{j}}\left(\dfrac{\partial^{2} y^\beta}{\partial
t^{2}}\right){\hat{g}_{\alpha\beta}} +{2}\dfrac{\partial}{\partial
x^{i}}\left(\dfrac{\partial y^{\alpha}}{\partial t}\right)
\dfrac{\partial y^{\beta}}{\partial x^{j}} \dfrac{d
\hat{g}_{\alpha\beta}}{d t}\vspace{2mm}\\
& & +{2}\dfrac{\partial y^{\alpha}}{\partial x^{i}}
\dfrac{\partial}{\partial x^{j}}\left(\dfrac{\partial
y^{\beta}}{\partial t}\right)\dfrac{d \hat{g}_{\alpha\beta}}{d
t}+{2}\dfrac{\partial}{\partial x^{i}}\left(\dfrac{\partial
y^{\alpha}}{\partial t}\right) \dfrac{\partial}{\partial
x^{j}}\left(\dfrac{\partial y^{\beta}}{\partial
t}\right)\hat{g}_{\alpha\beta}.
\end{array}\end{eqnarray}
On the other hand,
$$
\dfrac{d \hat{g}_{\alpha\beta}}{d t}(y(x,t),t) = \dfrac{\partial
\hat{g}_{\alpha\beta}}{\partial y^{\gamma}} \dfrac{\partial
y^{\gamma}}{\partial t}+\dfrac{\partial
\hat{g}_{\alpha\beta}}{\partial t},$$
$$\dfrac{d^{2}\hat{g}_{\alpha\beta}}{d t^{2}}(y(x,t),t) =
\dfrac{\partial^{2} \hat{g}_{\alpha\beta}}{\partial y^{\gamma}
\partial y^{\lambda}} \dfrac{\partial y^{\gamma}}{\partial
t} \dfrac{\partial y^{\lambda}}{\partial t}+{2}\dfrac{\partial^{2}
\hat{g}_{\alpha\beta}}{\partial y^{\gamma}
\partial t} \dfrac{\partial y^{\gamma}}{\partial
t}+\dfrac{\partial^{2} \hat{g}_{\alpha\beta}}{\partial
t^{2}}+\dfrac{\partial \hat{g}_{\alpha\beta}}{\partial y^{\gamma}}
\dfrac{\partial^{2} y^{\gamma}}{\partial t^{2}}$$ and
$$\dfrac{\partial^{2} \hat{g}_{\alpha\beta}}{\partial t^2}(y,t) =
-{2}\hat{R}_{\alpha\beta}(y,t).$$ It follows from (2.3) that
\begin{equation}\begin{array}{lll} \dfrac{\partial^{2} g_{ij}}{\partial t^{2}}(x,t)&=&
-{2}\hat{R}_{\alpha\beta}(y,t)\dfrac{\partial y^{\alpha}}{\partial
x^{i}} \dfrac{\partial y^{\beta}}{\partial
x^{j}}+\dfrac{\partial^{2} \hat{g}_{\alpha\beta}}{\partial
y^{\gamma}
\partial y^{\lambda}}\dfrac{\partial y^{\alpha}}{\partial
x^{i}} \dfrac{\partial y^{\beta}}{\partial x^{j}} \dfrac{\partial
y^{\gamma}}{\partial t} \dfrac{\partial y^{\lambda}}{\partial
t}\vspace{2mm}\\& &+ {2}\dfrac{\partial^{2}
\hat{g}_{\alpha\beta}}{\partial y^{\gamma}
\partial t} \dfrac{\partial y^{\alpha}}{\partial
x^{i}} \dfrac{\partial y^{\beta}}{\partial x^{j}} \dfrac{\partial
y^{\gamma}}{\partial t}+\dfrac{\partial}{\partial
x^{i}}\left(\hat{g}_{\alpha\beta}\dfrac{\partial y^{\beta}}{\partial
x^{j}}\dfrac{\partial^{2} y^{\alpha}}{\partial
t^{2}}\right)+\dfrac{\partial}{\partial
x^{j}}\left(\hat{g}_{\alpha\beta}\dfrac{\partial y^{\beta}}{\partial
x^{i}} \dfrac{\partial^{2} y^{\alpha}}{\partial
t^{2}}\right)\vspace{2mm}\\& &
+\left[\dfrac{\partial\hat{g}_{\alpha\beta}}{\partial y^{\gamma}}
\dfrac{\partial y^{\alpha}}{\partial x^{i}} \dfrac{\partial
y^{\beta}}{\partial x^{j}}-\dfrac{\partial}{\partial
x^{i}}\left(\dfrac{\partial y^{\beta}}{\partial x^{j}}\hat{g}_{\beta
\gamma}\right)-\dfrac{\partial}{\partial x^{j}}\left(\dfrac{\partial
y^{\beta}}{\partial x^{i}}\hat{g}_{\beta \gamma}\right)\right]
\dfrac{\partial^{2} y^{\gamma}}{\partial t^{2}}\vspace{2mm}\\& &
+{2}\dfrac{\partial}{\partial x^{i}}\left(\dfrac{\partial
y^\alpha}{\partial t}\right) \dfrac{\partial y^{\beta}}{\partial
x^{j}}\left(\dfrac{\partial \hat{g}_{\alpha\beta}}{\partial
t}+\dfrac{\partial \hat{g}_{\alpha\beta}}{\partial y^{\gamma}}
\dfrac{\partial y^{\gamma}}{\partial t}\right)+{2}\dfrac{\partial
y^{2}}{\partial x^{i}} \dfrac{\partial}{\partial
x^{j}}\left(\dfrac{\partial y^{\beta}}{\partial t}\right)
\left(\dfrac{\partial \hat{g}_{\alpha\beta}}{\partial y^{\gamma}}
\dfrac{\partial y^{\gamma}}{\partial t}+\dfrac{\partial
\hat{g}_{\alpha\beta}}{\partial t}\right)\vspace{2mm}\\& &
+{2}\hat{g}_{\alpha\beta}\dfrac{\partial}{\partial
x^{i}}\left(\dfrac{\partial y^{\alpha}}{\partial t}\right)
\dfrac{\partial}{\partial x^{j}}\left(\dfrac{\partial
y^{\beta}}{\partial t}\right).\end{array}\end{equation}

Let us choose the normal coordinates $\{x^{i}\}$ around a fixed
point $p\in{\mathscr{M}}$ such that $ \dfrac{\partial
g_{ij}}{\partial x^{k}}=0 $ at $p$. We next prove that, at
$p\in{\mathscr{M}}$,
\begin{equation}
\dfrac{\partial \hat{g}_{\alpha\beta}}{\partial y^{\gamma}}
\dfrac{\partial y^{\alpha}}{\partial x^{i}} \dfrac{\partial
y^{\beta}}{\partial x^{j}}-\dfrac{\partial}{\partial
x^{i}}\left(\dfrac{\partial y^{\beta}}{\partial x^{j}}\hat{g}_{\beta
\gamma}\right)-\dfrac{\partial}{\partial x^{j}}\left(\dfrac{\partial
y^{\beta}}{\partial x^{i}}\hat{g}_{\beta\gamma}\right)=0,\quad
\forall \;\; i,j,\gamma=1,\cdots,n.
\end{equation}

The left hand side of (2.5) is
\begin{eqnarray*}
& &\dfrac{\partial \hat{g}_{\alpha\beta}}{\partial y^{\gamma}}
\dfrac{\partial y^{\alpha}}{\partial x^{i}} \dfrac{\partial
y^{\beta}}{\partial x^{j}}-\dfrac{\partial}{\partial
x^{i}}\left(\dfrac{\partial y^{\beta}}{\partial x^{j}}\hat{g}_{\beta
\gamma}\right)- \dfrac{\partial}{\partial
x^{j}}\left(\dfrac{\partial y^{\beta}}{\partial
x^{i}}\hat{g}_{\alpha\beta}\right)\vspace{2mm}\\&=&
\dfrac{\partial}{\partial y^{\gamma}}\left(g_{mn}\dfrac{\partial
x^{m}}{\partial y^{\alpha}} \dfrac{\partial x^{n}}{\partial
y^{\beta}}\right) \dfrac{\partial y^{\alpha}}{\partial x^{i}}
\dfrac{\partial y^{\beta}}{\partial x^{j}}
-\dfrac{\partial}{\partial x^{i}}\left(g_{mj}\dfrac{\partial
x^{m}}{\partial y^{\gamma}}\right)-\dfrac{\partial}{\partial
x^{j}}\left(g_{mi}\dfrac{\partial x^{m}}{\partial
y^{\gamma}}\right)\\&=&g_{mn}\left[\dfrac{\partial}{\partial
y^{\gamma}}\left(\dfrac{\partial x^{m}}{\partial y^{\alpha}}\right)
\dfrac{\partial x^{n}}{\partial y^{\beta}} \dfrac{\partial
y^{\alpha}}{\partial x^{i}} \dfrac{\partial y^{\beta}}{\partial
x^{j}} + \dfrac{\partial x^{m}}{\partial y^{\alpha}}
\dfrac{\partial}{\partial y^{\gamma}}\left(\dfrac{\partial
x^{n}}{\partial y^{\beta}}\right) \dfrac{\partial
y^{\alpha}}{\partial x^{i}} \dfrac{\partial y^{\beta}}{\partial
x^{j}}\right]\vspace{2mm}\\& & -g_{mj}\dfrac{\partial}{\partial
x^{i}}\left(\dfrac{\partial x^{m}}{\partial
y^{\gamma}}\right)-g_{mi}\dfrac{\partial}{\partial
x^{j}}\left(\dfrac{\partial x^{m}}{\partial y^{\gamma}}\right)\\&=&
g_{mj}\dfrac{\partial}{\partial y^{\gamma}}\left(\dfrac{\partial
x^{m}}{\partial y^{\alpha}}\right) \dfrac{\partial
y^{\alpha}}{\partial
 x^{i}}+g_{mi}\dfrac{\partial}{\partial y^{\gamma}}\left(\dfrac{\partial
x^{m}}{\partial y^{\beta}}\right) \dfrac{\partial
y^{\beta}}{\partial
 x^{j}}-g_{mj}\dfrac{\partial}{\partial x^{i}}\left(\dfrac{\partial
x^{m}}{\partial y^{\gamma}}\right) -g_{mi}\dfrac{\partial}{\partial
x^{j}}\left(\dfrac{\partial x^{m}}{\partial
y^{\gamma}}\right)\vspace{2mm}\\&=& g_{mj}\dfrac{\partial}{\partial
y^{\alpha}}\left(\dfrac{\partial x^{m}}{\partial y^{\gamma}}\right)
\dfrac{\partial y^{\alpha}}{\partial
 x^{i}}+g_{mi}\dfrac{\partial}{\partial y^{\beta}}\left(\dfrac{\partial
x^{m}}{\partial y^{\gamma}}\right) \dfrac{\partial
y^{\beta}}{\partial
 x^{j}}
 -g_{mj}\dfrac{\partial}{\partial x^{i}}\left(\dfrac{\partial
x^{m}}{\partial y^{\gamma}}\right)-g_{mi}\dfrac{\partial}{\partial
x^{j}}\left(\dfrac{\partial x^{m}}{\partial
y^{\gamma}}\right)\\&=& 0.\end{eqnarray*} So (2.5) holds.

By (2.4) and (2.5), we have
\begin{eqnarray} \dfrac{\partial^{2} g_{ij}}{\partial t^{2}}(x,t)
&=& -{2}{R_{ij}}(x,t)+\dfrac{\partial}{\partial
x^{i}}\left(g_{mj}\dfrac{\partial x^{m}}{\partial y^{\alpha}}
\dfrac{\partial^{2} y^{\alpha}}{\partial
t^{2}}\right)+\dfrac{\partial}{\partial
x^{j}}\left(g_{mi}\dfrac{\partial y^{\alpha}}{\partial x^{m}}
\dfrac{\partial^{2} y^{\alpha}}{\partial t^{2}}\right)\nonumber
\vspace{2mm}\\& &
+\dfrac{\partial^{2}\hat{g}_{\alpha\beta}}{\partial y^{\gamma}
\partial y^{\lambda}} \dfrac{\partial y^{\alpha}}{\partial
x^{i}} \dfrac{\partial y^{\beta}}{\partial x^{j}} \dfrac{\partial
y^{\gamma}}{\partial t} \dfrac{\partial y^{\lambda}}{\partial
t}+{2}\dfrac{\partial \hat{g}_{\alpha\beta}}{\partial y^{\gamma}
\partial t} \dfrac{\partial y^{\gamma}}{\partial
t} \dfrac{\partial y^{\alpha}}{\partial x^{i}} \dfrac{\partial
y^{\beta}}{\partial x^{j}}\nonumber\vspace{2mm}\\& &+
2\dfrac{\partial }{{\partial x^i }}\left(\dfrac{{\partial y^\alpha
}}{{\partial t}}\right)\dfrac{{\partial y^\beta  }}{{\partial x^j
}}\left(\dfrac{{\partial \widehat g_{\alpha \beta } }}{{\partial t}}
+ \dfrac{{\partial \widehat g_{\alpha \beta } }}{{\partial y^\gamma
}}\dfrac{{\partial y^\gamma }}{{\partial t}}\right) +
2\dfrac{{\partial y^\alpha  }}{{\partial x^i }}\dfrac{\partial
}{{\partial x^j }}\left(\dfrac{{\partial y^\beta }}{{\partial
t}}\right)\left(\dfrac{{\partial \widehat g_{\alpha \beta }
}}{{\partial t}} + \dfrac{{\partial \widehat g_{\alpha \beta }
}}{{\partial y^\gamma }}\dfrac{{\partial
y^\gamma  }}{{\partial t}}\right) \nonumber\vspace{2mm} \\
&& +2\widehat g_{\alpha \beta } \dfrac{\partial }{{\partial x^i
}}\left(\dfrac{{\partial y^\alpha  }}{{\partial t}}\right)
\dfrac{\partial }{{\partial x^j }}\left(\dfrac{{\partial y^\beta
}}{{\partial t}}\right).
\end{eqnarray}
We define $y(x,t) = \varphi_t(x)$ by the following initial value
problem
\begin{equation}
\left\{ {\begin{array}{l}
   {\dfrac{{\partial ^2 y^\alpha  }}{{\partial t^2 }} =
   \dfrac{{\partial y^\alpha  }}{{\partial x^k }}g^{jl} (\Gamma _{jl}^k  - \overset{\circ\;\;\,}{{\Gamma}_{jl}^k} )},
   \vspace{2mm}  \\
   {y^\alpha  (x,0) = x^\alpha ,\quad \dfrac{\partial }{{\partial t}}y^\alpha  (x,0) = y_1^\alpha  (x)}  \\
\end{array}} \right. \end{equation}
and define the vector filed
$$V_i = g_{ik}g^{jl}\left(\Gamma _{jl}^k  -
\overset{\circ\;\;\,}{{\Gamma}_{jl}^k}\right),$$ where $\Gamma
_{jl}^k$ and $\overset{\circ\;\;\,}{{\Gamma}_{jl}^k}$ are the
connection coefficients corresponding to the metrics $g_{ij}(x,t)$
and $g_{ij}(x,0)$, respectively, $y_1^\alpha
(x)\;\;(\alpha=1,2,\cdots,n)$ are arbitrary $C^{\infty}$ smooth
functions on the manifold $\mathscr{M}$. We get the following
evolution equation for the pull-back metric
\begin{equation}\left\{\begin{array}{lll} \dfrac{\partial^{2}g_{ij}}{\partial
t^{2}}(x,t)&=&-2R_{ij}(x,t)+\nabla_{i}V_{j}+\nabla_{j}V_{i}+\dfrac{\partial^{2}\widehat{g}_{\alpha\beta}}{\partial
y^{\gamma}\partial y^{\lambda}} \dfrac{\partial y^{\alpha}}{\partial
x^{i}} \dfrac{\partial y^{\beta}}{\partial x^{j}} \dfrac{\partial
y^{\gamma}}{\partial t} \dfrac{\partial y^{\lambda}}{\partial
t}\vspace{2mm}\\
&&+2\dfrac{\partial^{2}\widehat{g}_{\alpha\beta}}{\partial
y^{\gamma}\partial t} \dfrac{\partial y^{\alpha}}{\partial x^{i}}
\dfrac{\partial y^{\beta}}{\partial x^{j}} \dfrac{\partial
y^{\gamma}}{\partial t}+2\dfrac{\partial}{\partial
x^{i}}\left(\dfrac{\partial y^{\alpha}}{\partial t}\right)
\dfrac{\partial y^{\beta}}{\partial
x^{j}}\left(\dfrac{\partial\widehat{g}_{\alpha\beta}}{\partial
t}+\dfrac{\partial\widehat{g}_{\alpha\beta}}{\partial y^{\gamma}}
\dfrac{\partial y^{\gamma}}{\partial t}\right)
\vspace{2mm}\\
&&+2\dfrac{\partial}{\partial x^{j}}\left(\dfrac{\partial
y^{\beta}}{\partial t}\right) \dfrac{\partial y^{\alpha}}{\partial
x^{i}}\left(\dfrac{\partial\widehat{g}_{\alpha\beta}}{\partial
t}+\dfrac{\partial\widehat{g}_{\alpha\beta}}{\partial y^{\gamma}}
\dfrac{\partial y^{\gamma}}{\partial
t}\right)\vspace{2mm}\\
&& +2\widehat{g}_{\alpha\beta}\dfrac{\partial}{\partial
x^{i}}\left(\dfrac{\partial y^{\alpha}}{\partial t}\right)
\dfrac{\partial}{\partial x^{j}}\left(\dfrac{\partial
y^{\beta}}{\partial t}\right)\vspace{2mm}\\
&\triangleq&-2R_{ij}(x,t)+\nabla_{i}V_{j}+\nabla_{j}V_{i}+F(Dy,D_{t}D_{x}y),
\vspace{2mm}\\
g_{ij}(x,0)&=&g^0_{ij}(x), \quad \quad\dfrac{\partial}{\partial
t}g_{ij}(x,0) \ = \ k^0_{ij}(x),
\end{array}\right.\end{equation} where
$$Dy=\left(\dfrac{\partial y^{\alpha}}{\partial t}, \dfrac{\partial
y^{\alpha}}{\partial x^{i}}\right), \quad
D_{t}D_{x}y=\left(\dfrac{\partial^{2}y^{\alpha}}{\partial
x^{i}\partial t}\right)\quad (\alpha, i=1,2,\cdots,n).$$

Let
\begin{equation*}
\widehat{\lambda}=\left(\dfrac{\partial y^{\alpha}}{\partial t},
\dfrac{\partial y^{\alpha}}{\partial x^{i}}, \dfrac{\partial^{2}
y^{\alpha}}{\partial x^{i}\partial t}\right)\quad (\alpha,
i=1,2,\cdots,n).
\end{equation*}
The nonlinear term $F=F(\widehat{\lambda})=F(Dy,D_{t}D_{x}y)$ in
(2.8) is smooth and
$F(\widehat{\lambda})=O(|\widehat{\lambda}|^{2})$ holds.

Since
\begin{equation*}
\Gamma^{k}_{ji}=\dfrac{\partial y^{\alpha}}{\partial x^{j}}
\dfrac{\partial y^{\beta}}{\partial x^{i}}\dfrac{\partial
x^{k}}{\partial
y^{\gamma}}\widehat{\Gamma}^{\gamma}_{\alpha\beta}+\dfrac{\partial
x^{k}}{\partial y^{\alpha}} \dfrac{\partial^{2}y^{\alpha}}{\partial
x^{j}\partial x^{i}},
\end{equation*}
the initial value problem (2.7) can be written as
\begin{equation}\left\{\begin{array}{l}
\dfrac{\partial^{2}y^{\alpha}}{\partial
t^{2}}=g^{jl}\left(\dfrac{\partial^{2}y^{\alpha}}{\partial
x^{j}\partial
x^{l}}-\overset{\circ\;\;\,}{{\Gamma}_{jl}^k}\dfrac{\partial
y^{\alpha}}{\partial
x^{k}}+\widehat{\Gamma}^{\alpha}_{\beta\gamma}\dfrac{\partial
y^{\beta}}{\partial x^{j}} \dfrac{\partial y^{\gamma}}{\partial
x^{i}}\right), \vspace{2mm}\\
y^{\alpha}(x,0)=x^{\alpha}, \qquad \dfrac{\partial}{\partial
t}y^{\alpha}(x,0) =  y^{\alpha}_{1}(x).
\end{array} \right.
\end{equation}
At the same time, in the normal coordinates $\{x^{i}\}$,
\begin{eqnarray*}
-2R_{ij}(x,t)+\nabla_{i}V_{j}+\nabla_{j}V_{i}
&=&\dfrac{\partial}{\partial x^{i}}\left\{g^{kl}\dfrac{\partial
g_{kl}}{\partial x^{j}}\right\}-\dfrac{\partial}{\partial
x^{k}}\left\{g^{kl}\left(\dfrac{\partial g_{jl}}{\partial
x^{i}}+\dfrac{\partial g_{il}}{\partial x^{j}}-\dfrac{\partial
g_{ij}}{\partial x^{l}}\right)\right\} \vspace{2mm}\\
&&+g_{jk}g^{pq}\dfrac{\partial}{\partial
x^{i}}\left\{\dfrac{1}{2}g^{kl}\left(\dfrac{\partial
g_{pl}}{\partial x^{q}}+\dfrac{\partial g_{ql}}{\partial
x^{p}}-\dfrac{\partial g_{pq}}{\partial
x^{l}}\right)\right\}\\&&+g_{ik}g^{pq}\dfrac{\partial}{\partial
x^{j}}\left\{\dfrac{1}{2}g^{kl}\left(\dfrac{\partial
g_{pl}}{\partial x^{q}}\dfrac{\partial g_{ql}}{\partial
x^{p}}-\dfrac{\partial g_{pq}}{\partial
x^{l}}\right)\right\}+({\rm
lower \ order \ terms})\vspace{2mm} \\
&=&g^{kl}\left\{\dfrac{\partial^{2}g_{kl}}{\partial x^{i}\partial
x^{j}}-\dfrac{\partial^{2}g_{jl}}{\partial x^{i}\partial
x^{k}}-\dfrac{\partial^{2}g_{il}}{\partial x^{j}\partial
x^{k}}+\dfrac{\partial^{2}g_{ij}}{\partial x^{k}\partial
x^{l}}\right\}\\&&+\dfrac{1}{2}g^{pq}\left\{\dfrac{\partial^{2}g_{pj}}{\partial
x^{i}\partial x^{q}}+\dfrac{\partial^{2}g_{qj}}{\partial
x^{i}\partial x^{q}}-\dfrac{\partial^{2} g_{pq}}{\partial
x^{i}\partial x^{j}}\right\} \vspace{2mm}\\
&&+\dfrac{1}{2}g^{pq}\left\{\dfrac{\partial^{2} g_{pi}}{\partial
x^{i}\partial x^{q}}+\dfrac{\partial^{2}g_{qi}}{\partial
x^{j}\partial x^{p}}-\dfrac{\partial^{2} g_{pq}}{\partial
x^{i}\partial x^{j}}\right\}+({\rm lower \ order \ terms}) \vspace{2mm}\\
&=&g^{kl}\dfrac{\partial^{2}g_{ij}}{\partial x^{k}\partial
x^{l}}+({\rm lower \ order \ terms}).
\end{eqnarray*}
Thereby, the initial value problem (2.8) can be written as
\begin{equation}\left\{ {\begin{array}{l}
\dfrac{\partial^{2}g_{ij}}{\partial
t^{2}}(x,t)=g^{kl}\dfrac{\partial^{2} g_{ij}}{\partial
x^{k}\partial
x^{l}}+F(Dy,D_{t}D_{x}y)+G(g,D_{x}g),\vspace{2mm}\\
g_{ij}(x,0)=\;g^0_{ij}(x), \qquad \dfrac{\partial}{\partial
t}g_{ij}(x,0) = k^0_{ij}(x),\end{array}} \right.
\end{equation}
where $g=(g_{ij})$, $D_{x}g=\left(\dfrac{\partial g_{ij}}{\partial
x^{k}}\right)$ ($i,j,k=1,2,\cdots,n)$. Let
\begin{equation*}
\widehat{\mu}=\left(g_{ij}, \dfrac{\partial g_{ij}}{\partial
x^{k}}\right) \quad (i,j,k=1,2,\cdots,n).
\end{equation*}
The nonlinear term $G=G(\widehat{\mu})=G(g,D_{x}g)$ in (2.10) is
smooth and quadratic with respect to $D_{x}g$.

We observe that both (2.9) and (2.10) are clearly strictly
hyperbolic systems. Since the equations (2.9) and (2.10) are
strictly hyperbolic and the manifold $\mathscr{M}$ is compact, it
follows from the standard theory of hyperbolic equations (see
\cite{jo}, \cite{ka}, \cite{k}) that the system united by (2.9) and
(2.10) has a unique smooth solution for a short time. Thus, we have
proved Theorem 1.1.

\section{Symmetrization of hyperbolic geometric flow --- second proof of Theorem 1.1}
In this section we reduce the hyperbolic geometric flow (2.1) to a
symmetric hyperbolic system. Then we use the theory of symmetric
hyperbolic system to give another proof of Theorem 1.1.

Let $\mathscr{M}$ be a compact $n$-dimensional manifold and
$g_{ij}(x,t)$ is a hyperbolic geometric flow on $\mathscr{M}$. We
denote the corresponding connection coefficients, the Riemannian
curvature tensor and the Ricci curvature tensor by $\Gamma
_{ij}^{k}$, $R_{ijkl}$ and $R_{ik}$, respectively.

We consider the space-time $\mathbb{R} \times \mathscr{M}$ with the
Lorentzian metric
\begin{equation}
ds^{2}=-dt^{2}+g_{ij}(x,t)dx^{i}dx^{j}
\end{equation}
and denote the corresponding connection coefficients, the
Riemannian curvature tensor and the Ricci curvature tensor by
$\widehat{\Gamma}_{\alpha \beta}^{\gamma}$, $\widehat R_{\alpha
\beta \gamma \lambda}$ and $\widehat R_{\alpha \beta}$,
respectively. Here and hereafter, the Greek indices run from $0$
to $n$, Latin indices from $1$ to $n$. The summation convention is
employed. We also denote $x^{0}=t$.

By direct computations,
\begin{equation}\begin{array}{l}
\widehat\Gamma_{ij}^{k}=\Gamma_{ij}^{k},\;\;\;
\widehat{\Gamma}_{ij}^{0}=\dfrac{1}{2} \dfrac{\partial
g_{ij}}{\partial t},\;\;\;
\widehat{\Gamma}_{0l}^{k}=\widehat{\Gamma}_{l0}^{k}=\dfrac{1}{2}g^{ki}\dfrac{\partial
g_{il}}{\partial t},\;\;\;\widehat{\Gamma}
_{0i}^{0}=0,\;\;\;\widehat{\Gamma} _{00}^{i}=0,\;\;\;
\widehat{\Gamma}_{00}^{0}=0,\vspace{3mm}\\
\widehat R_{ijl}^{k}=\dfrac{\partial \Gamma_{jl}^{k}}{\partial
x^{i}}-\dfrac{\partial \Gamma_{il}^{k}}{\partial
x^{j}}+\widehat{\Gamma}_{i\alpha}^{k}\widehat{\Gamma}_{jl}^{\alpha}-\widehat{\Gamma}
_{j\alpha}^{k}\widehat{\Gamma}_{il}^{\alpha} =
R_{ijl}^{k}+\dfrac{1}{4}g^{km}\dfrac{\partial g_{mi}}{\partial
t}\dfrac{\partial g_{jl}}{\partial t} -
\dfrac{1}{4}g^{km}\dfrac{\partial g_{mj}}{\partial
t}\dfrac{\partial g_{il}}{\partial t},\vspace{3mm}\\ \widehat
R_{0k0}^{l}=\dfrac{1}{2}\dfrac{\partial g^{li}}{\partial t}
\dfrac{\partial g_{ik}}{\partial t}+
\dfrac{1}{2}g^{li}\dfrac{\partial^{2}g_{ik}}{\partial t^{2}}+
\dfrac{1}{4}g^{lm}g^{pn}\dfrac{\partial g_{mp}}{\partial
t}\dfrac{\partial g_{nk}}{\partial t}.
\end{array}\end{equation}
Then,
\begin{eqnarray}
\nonumber \widehat R_{ik}&=& g^{\alpha \beta}\widehat R_{i\alpha k
\beta}=g^{\alpha \beta}g_{kl}\widehat
R_{i\alpha\beta}^{l}=-g_{kl}\widehat
R_{i00}^{l}+g^{pq}g_{kl}\widehat
               R_{ipq}^{l}\\\nonumber
               &=&-g_{kl}\left[-\dfrac{1}{2}\dfrac{\partial g^{lp}}{\partial t}\dfrac{\partial g_{ip}}
               {\partial t}-\dfrac{1}{2}g^{lp}
               \dfrac{\partial^{2}g_{ip}}{\partial t^{2}}-\dfrac{1}{4}g^{lm}g^{pn}\dfrac{\partial g_{mp}}{\partial t}
               \dfrac{\partial g_{ni}}{\partial t}\right]\\\nonumber
               &&+ g^{pq}g_{kl}\left[R_{ipq}^{l}+\dfrac{1}{4}g^{lm}
               \dfrac{\partial g_{mi}}{\partial t}\dfrac{\partial g_{pq}}{\partial t}-
               \dfrac{1}{4}g^{lm}\dfrac{\partial g_{mp}}{\partial t} \dfrac{\partial g_{iq}}{\partial
               t}\right]\\\nonumber
               &=&\dfrac{1}{2}\dfrac{\partial^{2} g_{ik}}{\partial
               t^{2}}+R_{ik}+\dfrac{1}{2}g_{kl}\dfrac{\partial g^{lp}}{\partial
               t}\dfrac{\partial g_{ip}}{\partial
               t}+\dfrac{1}{4}g^{pn}\dfrac{\partial g_{kp}}{\partial
               t}\dfrac{\partial g_{ni}}{\partial
               t}+\dfrac{1}{4}g^{pq}\dfrac{\partial g_{ki}}{\partial
               t}\dfrac{\partial g_{pq}}{\partial
               t}-\dfrac{1}{4}g^{pq}\dfrac{\partial g_{kp}}{\partial
               t}\dfrac{\partial g_{iq}}{\partial t}\\
               &=&\dfrac{1}{2}\dfrac{\partial^{2}g_{ik}}{\partial
               t^{2}}+R_{ik}+\dfrac{1}{4}g^{pq}\dfrac{\partial g_{ik}}{\partial
               t}\dfrac{\partial g_{pq}}{\partial t}-\dfrac{1}{2}g^{pq}\dfrac{\partial g_{ip}}{\partial
               t}\dfrac{\partial g_{kq}}{\partial t}.
          \end{eqnarray}
A direct computation gives (see, e.g., Fock \cite{fo}, p.423;
Fisher and Marsden \cite{fm}, p.22)
$$ \widehat R_{ij}=\widehat
R_{ij}^{(h)}+\dfrac{1}{2}\left(g_{i\alpha}\dfrac{\partial \widehat
\Gamma^{\alpha}}{\partial x^{j}}+g_{j\alpha}\dfrac{\partial
\widehat \Gamma^{\alpha}}{\partial x^{i}}\right)=\widehat
R_{ij}^{(h)}+\dfrac{1}{2}\left(g_{ik}\dfrac{\partial
\Gamma^{k}}{\partial x^{j}}+g_{jk}\dfrac{\partial
\Gamma^{k}}{\partial x^{i}}\right),
$$
where
$$\widehat{\Gamma}^{\alpha}=g^{\beta\gamma}\widehat{\Gamma}^{\alpha}_{\beta\gamma},\;\;\mbox{i.e.,}\;\;
\widehat{\Gamma}^{0}=g^{kl}\widehat{\Gamma}^{0}_{kl}=\frac{1}{2}g^{kl}\dfrac{\partial
g_{kl}}{\partial
t},\;\;\widehat{\Gamma}^{i}=g^{kl}\widehat{\Gamma}^{i}_{kl}=g^{kl}\Gamma^{i}_{kl}\stackrel{\triangle}{=}\Gamma^i,$$
$$
\widehat R_{ij}^{(h)}=-\dfrac {1}{2}g^{\alpha
\beta}\dfrac{\partial^{2}g_{ij}}{\partial x^{\alpha}\partial
x^{\beta}}+\widehat H_{ij}\left(g_{\alpha \beta},\dfrac{\partial
g_{\alpha \beta}}{\partial
x^{\lambda}}\right)=\dfrac{1}{2}\dfrac{\partial^{2}g_{ij}}{\partial
t^{2}}-\dfrac{1}{2}g^{kl}\dfrac{\partial^{2}g_{ij}}{\partial
x^{k}\partial x^{l}}+\widehat H_{ij}\left(g_{\alpha
\beta},\dfrac{\partial g_{\alpha\beta}}{\partial
x^{\lambda}}\right)
$$
and
$$
\widehat H_{ij}\left(g_{\alpha \beta},\dfrac{\partial g_{\alpha
\beta}}{\partial x^{\lambda}}\right)=g^{\alpha
\beta}g_{e\sigma}\widehat \Gamma_{i\beta}^{e}\widehat \Gamma
_{j\alpha}^{\sigma}+\dfrac{1}{2}\left(\dfrac{\partial
g_{ij}}{\partial x^{\alpha}}\widehat
\Gamma^{\alpha}+g_{j\lambda}\widehat \Gamma _{\alpha \beta}
^{\lambda}g^{\alpha \gamma}g^{\beta g}\dfrac{\partial g_{\gamma
g}}{\partial x^{i}}+g_{i\lambda}\widehat \Gamma_{\alpha
\beta}^{\lambda}g^{\alpha \gamma}g^{\beta g}\dfrac{\partial
g_{\gamma g}}{\partial x^{j}}\right).
$$
Similar to (3.3), we have
$$
\widehat R_{ij}^{(h)}=\dfrac{1}{2}\dfrac{\partial
^{2}g_{ij}}{\partial t^{2}}-\dfrac{1}{2}g^{kl}\dfrac{\partial
^{2}g_{ij}}{\partial x^{k}\partial
x^{l}}+H_{ij}\left(g_{kl},\dfrac{\partial g_{kl}}{\partial
x^{p}},\dfrac{\partial g_{kl}}{\partial t}\right)
$$
and
\begin{eqnarray}
H_{ij}\left(g_{kl},\dfrac{\partial g_{kl}}{\partial
x^{p}},\dfrac{\partial g_{kl}}{\partial t}\right)\nonumber &
\stackrel {\Delta}{=}& \widehat H_{ij}\left(g_{\alpha
\beta},\dfrac{\partial g_{\alpha \beta}}{\partial
x^{\lambda}}\right)\\\nonumber &=&-g_{kl}\widehat
\Gamma_{i0}^{k}\widehat\Gamma
_{j0}^{l}+(-1)g^{kl}\widehat\Gamma_{ik}^{0}\widehat\Gamma_{jl}^{0}+g^{kl}g_{pq}
\Gamma _{ik}^{p}\Gamma_{jl}^{q}\\\nonumber
&&+\dfrac{1}{2}\left(\dfrac{\partial g_{ij}}{\partial
t}\dfrac{1}{2}g^{kl}\dfrac{\partial g_{kl}}{\partial
t}+\dfrac{\partial g_{ij}}{\partial
x_k}\Gamma_{pq}^kg^{pq}+g_{jk}\Gamma _{rs}^{k}\dfrac{\partial
g_{pq}}{\partial x^{i}}g^{pr}g^{qs}+g_{ik}\Gamma
_{rs}^{k}g^{pr}g^{qs}\dfrac{\partial g_{pq}}{\partial
x^{j}}\right)\\\nonumber
&=&-g_{kl}\dfrac{1}{2}g^{kp}\dfrac{\partial g_{ip}}{\partial
t}\dfrac{1}{2}g^{lq}\dfrac{\partial g_{jq}}{\partial
t}-g^{kl}\dfrac{1}{2}\dfrac{\partial g_{ik}}{\partial
t}\dfrac{1}{2}\dfrac{\partial g_{jl}}{\partial
t}+g^{kl}g_{pq}\Gamma _{ik}^{p}\Gamma_{jl}^{q}\\\nonumber
 &&+
\dfrac{1}{4}g^{kl}\dfrac{\partial g_{kl}}{\partial
t}\dfrac{\partial g_{ij}}{\partial t}+\dfrac{1}{2}\dfrac{\partial
g_{ij}}{\partial x_k}\Gamma_{pq}^kg^{pq}\\\nonumber
&&+\dfrac{1}{2}\left(g_{ik}\Gamma
_{rs}^{k}g^{pr}g^{qs}\dfrac{\partial g_{pq}}{\partial
x^{j}}+g_{jk}\Gamma _{rs}^{k}g^{pr}g^{qs}\dfrac{\partial
g_{pq}}{\partial
x^{i}}\right)\\\nonumber&=&-\dfrac{1}{2}g^{pq}\dfrac{\partial
g_{ip}}{\partial t}\dfrac{\partial g_{jq}}{\partial t}+
\dfrac{1}{4}g^{kl}\dfrac{\partial g_{kl}}{\partial
t}\dfrac{\partial g_{ij}}{\partial t}\\\nonumber&
&+g^{kl}g_{pq}\Gamma
_{ik}^{p}\Gamma_{jl}^{q}+\dfrac{1}{2}\dfrac{\partial
g_{ij}}{\partial x_k}\Gamma_{pq}^kg^{pq}\\\nonumber
&&+\dfrac{1}{2}\left(g_{ik}\Gamma
_{rs}^{k}g^{pr}g^{qs}\dfrac{\partial g_{pq}}{\partial
x^{j}}+g_{jk}\Gamma _{rs}^{k}g^{pr}g^{qs}\dfrac{\partial
g_{pq}}{\partial x^{i}}\right).
\end{eqnarray}
It follows from (3.3) that
\begin{eqnarray*}
\dfrac{1}{2}\dfrac{\partial^{2} g_{ij}}{\partial
t^{2}}+R_{ij}&=&\dfrac{1}{2}\dfrac{\partial^{2}g_{ij}}{\partial
t^{2}}-\dfrac{1}{2}g^{kl}\dfrac{\partial^{2}g_{ij}}{\partial x^{k}
\partial x^{l}}+\dfrac{1}{2}\left(g_{ik}\dfrac{\partial \Gamma^{k}}{\partial
x^{j}}+g_{jk}\dfrac{\partial \Gamma^{k}}{\partial
x^{i}}\right)\\\nonumber & &+H_{ij}\left(g_{kl},\dfrac{\partial
g_{kl}}{\partial x^{p}},\dfrac{\partial g_{kl}}{\partial t}\right)
+\dfrac{1}{2}g^{pq}\dfrac{\partial g_{ip}}{\partial
t}\dfrac{\partial g_{jq}}{\partial
t}-\dfrac{1}{4}g^{pq}\dfrac{\partial g_{ij}}{\partial
t}\dfrac{\partial g_{pq}}{\partial t},
\end{eqnarray*}
i.e.,
\begin{eqnarray}
\nonumber\dfrac{\partial^{2} g_{ij}}{\partial
t^{2}}+2R_{ij}&=&\dfrac{\partial ^{2}g_{ij}}{\partial
t^{2}}-g^{kl}\dfrac{\partial^{2}g_{ij}}{\partial x^{k}\partial
x^{l}}+\left(g_{ik}\dfrac{\partial \Gamma^{k}}{\partial
x^{j}}+g_{jk}\dfrac{\partial \Gamma^{k}}{\partial
x^{i}}\right)\\\nonumber & &+2H_{ij}\left(g_{kl},\dfrac{\partial
g_{kl}}{\partial x^{p}},\dfrac{\partial g_{kl}}{\partial
t}\right)+g^{pq}\dfrac{\partial g_{ip}}{\partial t}\dfrac{\partial
g_{jq}}{\partial t}-\dfrac{1}{2}g^{pq}\dfrac{\partial
g_{ij}}{\partial t}\dfrac{\partial g_{pq}}{\partial
t}\\\nonumber&=&\dfrac{\partial ^{2}g_{ij}}{\partial
t^{2}}-g^{kl}\dfrac{\partial^{2}g_{ij}}{\partial x^{k}\partial
x^{l}}+\left(g_{ik}\dfrac{\partial \Gamma^{k}}{\partial
x^{j}}+g_{jk}\dfrac{\partial \Gamma^{k}}{\partial
x^{i}}\right)\\\nonumber & &+2g^{kl}g_{pq}\Gamma
_{ik}^{p}\Gamma_{jl}^{q}+\dfrac{\partial g_{ij}}{\partial
x_k}\Gamma_{pq}^kg^{pq}\\ &&+\left(g_{ik}\Gamma
_{rs}^{k}g^{pr}g^{qs}\dfrac{\partial g_{pq}}{\partial
x^{j}}+g_{jk}\Gamma _{rs}^{k}g^{pr}g^{qs}\dfrac{\partial
g_{pq}}{\partial x^{i}}\right).
\end{eqnarray}

Similar to the harmonic coordinates in the space-time (see \cite
{fm}), here we make use of a new kind of coordinates on the manifold
$\mathscr{M}$ defined by
\begin{equation}
\qquad \Gamma^{i}\stackrel {\Delta}{=}g^{kl}\Gamma_{kl}^{i}=0.
\end{equation}
Such new coordinates are called {\it elliptic coordinates} on
$\mathscr{M}$.

\begin{Lemma} Let $g_{ij}$ be a $C^{\infty}$ Riemannian metric on the manifold
$\mathscr{M}$. There is a $C^{\infty}$ local coordinates
transformation $\phi:
\;\;\mathscr{M}\longrightarrow\mathscr{M},\;\;x\longrightarrow\overline{x}$
around a fixed point $p\in\mathscr{M}$ such that the transformed
metric $\overline{g}_{ij}$ is a $C^{\infty}$ Riemannian metric with
$\overline{\Gamma}^{k}(\overline{x})=0$ for all $\overline{x}$ in a
neighborhood around $p\in\mathscr{M}$ and any
$k\in\{1,2,\cdots,n\}$.
\end{Lemma}

\noindent {\bf Proof.} Consider the elliptic equation for the
scalar $\psi$,
$$\triangle\psi\stackrel{\triangle}{=}-g^{kl}\dfrac{\partial^2\psi}{\partial x_k\partial
x_l}+g^{kl}\Gamma^{j}_{kl}\dfrac{\partial \psi}{\partial x_j}=0.$$
The coefficients are $C^{\infty}$. Let $\overline{x}^i(x)$ be a
solution with the condition
$\overline{x}^i(p)=x^i(p),\;\;\dfrac{\partial\overline{x}^i}{\partial
x_j}(p)=\delta^i_j$. Then $\overline{x}^i$ is a $C^{\infty}$ local
coordinates transformation $\phi:
\;\;\mathscr{M}\longrightarrow\mathscr{M}$ around $p$ and the
transformed metric $\overline{g}_{ij}$ is a $C^{\infty}$ Riemannian
metric.

Now the equation $\triangle\psi=0$ is a tensorial (scalar)
equation. In the barred coordinates, it becomes
$$0=\overline{\triangle}\overline{x}^k=-\overline{g}^{ij}
\dfrac{\partial^2\overline{x}^k}{\partial \overline{x}_i\partial
\overline{x}_j}+\overline{\Gamma}^{j}\dfrac{\partial
\overline{x}^k}{\partial \overline{x}_j}=\overline{\Gamma}^k.$$
Therefore, $\overline{g}_{ij}$ satisfies the elliptic condition
(3.5). The proof of Lemma 3.1 is complete.  $\quad\P$

\vskip 4mm

By Lemma 3.1, we can choose the elliptic coordinates around a
fixed point $p\in\mathscr{M}$ and for a fixed time
$t\in\mathbb{R}^+$. After throwing off the bar sign, the geometric
hyperbolic flow (3.1) can be written as
\begin{equation}
\dfrac{\partial ^{2}g_{ij}}{\partial
t^{2}}=g^{kl}\dfrac{\partial^{2}g_{ij}}{\partial x^{k}\partial
x^{l}}+\widetilde{H}_{ij}\left(g_{kl},\dfrac{\partial
g_{kl}}{\partial x^{p}}\right),
\end{equation}
where
\begin{equation}
\widetilde{H}_{ij}\left(g_{kl},\dfrac{\partial g_{kl}}{\partial
x^{p}}\right)=-2g^{kl}g_{pq}\Gamma^p_{ik}\Gamma^q_{jl}-\left(g_{ik}\Gamma
_{rs}^{k}g^{pr}g^{qs}\dfrac{\partial g_{pq}}{\partial
x^{j}}+g_{jk}\Gamma _{rs}^{k}g^{pr}g^{qs}\dfrac{\partial
g_{pq}}{\partial x^{i}}\right)
\end{equation}
are homogenous quadratic with respect to $\dfrac{\partial
g_{kl}}{\partial x^{p}}$ and rational with respect to $g_{kl}$
with non-zero denominator $\rm det (g)\neq0$. By introducing the
new unknowns $g_{ij},\ h_{ij}=\dfrac{\partial g_{ij}}{\partial
t},\ g_{ij,k}=\dfrac{\partial g_{ij}}{\partial x^k}$, the system
(3.6) can be reduced to a system of partial differential equations
of first order. We now consider such a quasi-linear (symmetric
hyperbolic) system with $\dfrac {1}{2}n(n+1)(n+2)$ PDEs of first
order
\begin{equation}
\begin{cases}
\dfrac{\partial g_{ij}}{\partial t}=h_{ij}, \vspace{2mm}\\
g^{kl}\dfrac{\partial g_{ij,k}}{\partial t}=g^{kl}\dfrac{\partial
h_{ij}}{\partial x^{k}},\vspace{2mm}\\
\dfrac{\partial h_{ij}}{\partial t}=g^{kl}\dfrac{\partial
g_{ij,k}}{\partial x^{l}}+\widetilde{H}_{ij}.
\end{cases}
\end{equation}
In the $C^2$ class, the system (3.8) is equivalent to (3.6).

Let $u=(g_{ij},g_{ij,k},h_{ij})^T$ be the
$\dfrac{1}{2}n(n+1)(n+2)$-dimensional unknown vector function. The
coefficient matrices $A^0,A^j,B$ are given by
$$
A^0(u)=A^0(g_{ij},g_{ij,k},h_{ij})=\left(\begin{array}{cccccc}I & 0 & 0 & \cdots & 0 & 0\\
 0 & g^{11}I & g^{12}I & \cdots & g^{1n}I & 0\\
 0 & g^{21}I & g^{22}I & \cdots & g^{2n}I & 0\\
 \vdots & \cdots\\\texttt{}
 0 & g^{n1}I & g^{n2}I & \cdots & g^{nn}I & 0\\
 0 & 0 & 0 & \cdots & 0 & I\end{array}\right),
$$
$$
A^{j}(u)=A^{j}(g_{kl},g_{kl,p},h_{kl})=\left(\begin{array}{cccccc}0
& 0 & 0 & \cdots & 0 & 0\\
0 & 0 & 0 & \cdots & 0 & g^{j1}I\\
0 & 0 & 0 & \cdots & 0 & g^{j2}I\\
\cdots & \cdots &\cdots\\
0 & 0 & 0 & \cdots & 0 & g^{jn}I\\
0 & g^{1j}I & g^{2j}I & \cdots & g^{nj}I & 0\end{array}\right),
$$
where 0 is the
$\left(\dfrac{1}{2}n(n+1)\right)\times\left(\dfrac{1}{2}n(n+1)\right)$
zero matrix, $I$ is the
$\left(\dfrac{1}{2}n(n+1)\right)\times\left(\dfrac{1}{2}n(n+1)\right)$
 identity matrix,
$$B(u)=B(g_{ij},g_{ij,p},h_{ij})=\left(\begin{array}{c}h_{ij}\\
0\\
 \widetilde{H}_{ij}\\\end{array}\right),$$
in which 0 is the $\dfrac{1}{2}n^2(n+1)$-dimensional zero vector.

We observe that the symmetric hyperbolic system
\begin{equation}
A^0(u)\dfrac{\partial u}{\partial t}=A^j(u)\dfrac{\partial
u}{\partial x^j}+B(u)
\end{equation}
is nothing but the system (3.8). So far, we have reduced the
equation of the hyperbolic geometric flow (3.1) to the symmetric
hyperbolic system (3.9), which are equivalent to each other in the
$C^2$ class. Then, by the theory of the symmetric hyperbolic
system, the smooth solution to (3.1) exists uniquely for a short
time (see \cite{fm}). Thus, the proof of Theorem 1.1 is completed.
\begin{Remark} The elliptic coordinates can also be used to prove
the short-time existence for the Ricci flow.
\end{Remark}

More generally, motivated by general Einstein equations and the
rich theory of hyperbolic equations, we may also consider the
following field equations with the energy-momentum tensor $T_{ij}$
under certain conditions:
\begin{equation}
\alpha_{ij}\frac{\partial^{2}g_{ij}}{\partial
t^{2}}+2R_{ij}+\mathscr{F}_{ij}\left(g,\frac{\partial g}{\partial
t}\right)=\kappa T_{ij},
\end{equation}
where $\alpha_{ij}$ are certain smooth functions on $\mathscr{M}$
which may also depend on $t$, $\mathscr{F}_{ij}$ are some given
smooth functions of the Riemannian metric $g$ and its first order
derivative with respect to $t$, and $\kappa$ is a parameter.
Similar results can be obtained.

\section{Nonlinear stability for hyperbolic geometric flow}

In this section we investigate the nonlinear stability of the
hyperbolic geometric flow defined on the Euclidean space with the
dimension larger than 4.

We now state the definition of nonlinear stability of the hyperbolic
geometric flow (1.1). Let $\mathscr{M}$ be an $n$-dimensional
complete Riemannian manifold. Given symmetric tensors $g_{ij}^{0}$
and $g_{ij}^{1}$ on $\mathscr{M}$, we consider the following initial
value problem
\begin{equation}
\begin{cases}
\dfrac{\partial^{2}}{\partial t^{2}}g_{ij}(x,t)=-2R_{ij}(x,t),\\
g_{ij}(x,0)=\overline{g}_{ij}(x)+\varepsilon
g_{ij}^{0}(x),\quad\dfrac{\partial g_{ij}}{\partial
t}(x,0)=\varepsilon g_{ij}^1(x),
\end{cases}
\end{equation}
where $\varepsilon >0$ is a small parameter.

\begin{Definition}
The Ricci flat Riemannian metric $\overline g _{ij}(x)$ possesses
the (locally) nonlinear stability with respect to
$(g_{ij}^{0},g_{ij}^{1})$, if there exists a positive constant
$\varepsilon_{0}=\varepsilon_{0}(g_{ij}^{0},g_{ij}^{1})$ such
that, for any $\varepsilon \in (0,\varepsilon_{0}]$, the initial
value problem (4.1) has a unique (local) smooth solution
$g_{ij}(x,t)$;

$\overline g_{ij}(x)$ is said to be (locally) nonlinearly stable, if
it possesses the (locally) nonlinear stability with respect to
arbitrary symmetric tensors $g_{ij}^{0}(x)$ and $g_{ij}^{1}(x)$ with
compact support.
\end{Definition}

In what follows, we consider the nonlinear stability of the flat
metric of the Euclidean space $\mathbb{R}^{n}$ with the dimension
$n\geq 5$. We have
\begin{Theorem}
The flat metric $g_{ij}=\delta_{ij}$ of the Euclidean space
$\mathbb{R}^{n}$ with $n\geq 5$ is nonlinearly stable.
\end{Theorem}
\begin{Remark}
Theorem 4.1 gives the nonlinear stability of the hyperbolic
geometric flow on the Euclidean space with dimension larger than 4.
The situation for the 3-, 4-dimensional Euclidean spaces is very
different, and will be studied in the sequel by using null
conditions. This is similar to the proofs of the Poincar$\acute{e}$
conjecture in topology: the proofs for the three and four
dimensional case and $n\ge 5$ dimensional case are very different
(see, for example, \cite{h} and \cite{cz}). This motivates us to
understand the possibility of using hyperbolic geometric flow to
understand Poincar$\acute{e}$ conjecture.
\end{Remark}
\noindent{\bf Proof of Theorem 4.1.} \ \ Define a 2-tensor $h$ in
the following way $$g_{ij}(x,t)=\delta_{ij}+ h_{ij}(x,t).$$ Let
$\delta^{ij}$ be the inverse of $\delta_{ij}$. Then, for small $h$
\begin{equation*}
H^{ij}\stackrel{\triangle}{=} g^{ij}-\delta^{ij}=-h^{ij}+
O^{ij}(h^{2}),\end{equation*} where
$h^{ij}=\delta^{ik}\delta^{jl}h_{kl}$ and $O^{ij}(h^{2})$ vanishes
to second order at $h=0$.

As in Section 3, we choose the above elliptic coordinates
$\{x^{i}\}$ around the origin in $\mathbb{R}^{n}$. Then the initial
value problem (4.1) can be written as
\begin{equation}
\begin{cases}
\dfrac{\partial^{2}}{\partial
t^{2}}h_{ij}(x,t)=\left(\delta^{kl}+H^{kl}\right)
\dfrac{\partial^{2}h_{ij}}{\partial x^{k}\partial x^{l}}+\tilde
{H}_{ij}
\left(\delta_{kl}+h_{kl}(x,t),\dfrac{\partial h_{kl}(x,t)}{\partial x^{p}}\right),\\
h_{ij}(x,0)=\varepsilon g_{ij}^{0}(x), \quad
\dfrac{\partial}{\partial t}h_{ij}(x,0)=\varepsilon g_{ij}^{1}(x),
\end{cases}
\end{equation}
where
\begin{eqnarray*}
&&\tilde H_{ij} \left(\delta_{kl}+h_{kl}(x,t),\dfrac{\partial
h_{kl}(x,t)}{\partial x^{p}}\right)=\tilde
H_{ij}\left(g_{kl}(x,t),\dfrac{\partial g_{kl}(x,t)}{\partial
x^{p}}\right)\\&=& -\left[2g^{kl}g_{pq}\Gamma
_{ik}^{p}\Gamma^{q}_{jl}+ g_{ik}\Gamma _{rs}^{k}g^{pr}g^{qs}
\dfrac{\partial g_{pq}}{\partial x^{j}}+g_{jk}\Gamma
_{rs}^{k}g^{pr}g^{qs}\dfrac{\partial
g_{pq}}{\partial x^{i}}\right] \\
&=& -2\left(\delta^{kl}+H^{ij}\right)\left(\delta
_{pq}+h_{pq}\right)\dfrac{1}{4}\left(\delta
^{pa}+H^{pa}\right)\left(\delta^{qb}+H^{qb}\right)\left(\dfrac{\partial
h_{ai}}{\partial x^{k}}+\dfrac{\partial h_{ak}}{\partial
x^{i}}-\dfrac{\partial h_{ik}}{\partial
x^{a}}\right)\left(\dfrac{\partial h_{bj}}{\partial
x^{l}}+\dfrac{\partial h_{bl}}{\partial x^{j}}-\dfrac{\partial
h_{jl}}{\partial x^{b}}\right)\\
&&-\left(\delta_{ik}+H_{ik}\right)\left(\delta^{pr}+H^{pr}\right)\left(\delta^{qs}+H^{qs}\right)
\dfrac{1}{2}\left(\delta^{ka}+H^{ka}\right)
\left(\dfrac{\partial h_{ar}}{\partial x^{s}}+\dfrac{\partial
h_{as}}{\partial x^{r}}-\dfrac{\partial h_{rs}}{\partial
x^{a}}\right)\dfrac{\partial h_{pq}}{\partial x^{j}} \\
&&-\left(\delta_{jk}+H_{jk}\right)\left(\delta^{pr}+H^{pr}\right)\left(\delta^{qs}+H^{qs}\right)
\dfrac{1}{2}\left(\delta^{ka}+H^{ka}\right)
\left(\dfrac{\partial h_{ar}}{\partial x^{s}}+\dfrac{\partial
h_{as}}{\partial x^{r}}-\dfrac{\partial h_{rs}}{\partial
x^{a}}\right)\dfrac{\partial h_{pq}}{\partial x^{i}} \\
&=&-2\cdot\dfrac{1}{4}\cdot\delta^{kl}\delta_{pq}\delta^{pa}\delta^{qb}\left(\dfrac{\partial
h_{ai}}{\partial x^{k}}+\dfrac{\partial h_{ak}}{\partial
x^{i}}-\dfrac{\partial h_{ik}}{\partial
x^{a}}\right)\left(\dfrac{\partial h_{bj}}{\partial
x^{l}}+\dfrac{\partial h_{bl}}{\partial
x^{j}}-\dfrac{\partial h_{jl}}{\partial x^{b}}\right)\\
&&-\dfrac{1}{2}\delta_{ik}\delta^{pr}\delta^{qs}\delta^{ka}\left(\dfrac{\partial
h_{ar}}{\partial x^{s}}+\dfrac{\partial h_{as}}{\partial
x^{r}}-\dfrac{\partial h_{rs}}{\partial
x^{a}}\right)\dfrac{\partial
h_{pq}}{\partial x^{j}}\\
&&-\dfrac{1}{2}\delta_{jk}\delta^{pr}\delta^{qs}\delta^{ka}\left(\dfrac{\partial
h_{ar}}{\partial x^{s}}+\dfrac{\partial h_{as}}{\partial
x^{r}}-\dfrac{\partial h_{rs}}{\partial x^{a}}\right)\dfrac{\partial
h_{pq}}{\partial x^{i}}+O\left(|h_{kl}|+\left|\dfrac{\partial
h_{kl}}{\partial x^{p}}\right|\right)^{3}\\&=&
-\dfrac{1}{2}\delta^{kl}\delta^{ab}\left(\dfrac{\partial
h_{ai}}{\partial x^{k}}+\dfrac{\partial h_{ak}}{\partial
x^{i}}-\dfrac{\partial h_{ik}}{\partial
x^{a}}\right)\left(\dfrac{\partial h_{bj}}{\partial
x^{l}}+\dfrac{\partial h_{bl}}{\partial
x^{j}}-\dfrac{\partial h_{jl}}{\partial x^{b}}\right)\\
&&-\dfrac{1}{2}\delta^{pr}\delta^{qs}\left(\dfrac{\partial
h_{ir}}{\partial x^{s}}+\dfrac{\partial h_{is}}{\partial
x^{r}}-\dfrac{\partial h_{rs}}{\partial
x^{i}}\right)\dfrac{\partial
h_{pq}}{\partial x^{j}}\\
&&-\dfrac{1}{2}\delta^{qs}\delta^{pr}\left(\dfrac{\partial
h_{jr}}{\partial x^{s}}+\dfrac{\partial h_{js}}{\partial
x^{r}}-\dfrac{\partial h_{rs}}{\partial x^{j}}\right)\dfrac{\partial
h_{pq}}{\partial x^{i}}+O\left(|h_{kl}|+\left|\dfrac{\partial
h_{kl}}{\partial x^{p}}\right|\right)^{3}.
\end{eqnarray*}
Furthermore, (4.2) can be written as
\begin{equation}
\begin{cases}
\dfrac{\partial^{2}}{\partial t^{2}}h_{ij}(x,t)={\displaystyle
\sum_{k=1}^{n}}\dfrac{\partial^{2}h_{ij}}{\partial x^{k}\partial
x^{k}}+ \bar{H}_{ij}\left(h_{kl},\dfrac{\partial
h_{kl}}{\partial x^{p}},\dfrac{\partial^{2}h_{kl}}{\partial x^{p}\partial x^{q}}\right),\vspace{2mm}\\
h_{ij}(x,0)=\varepsilon
g_{ij}^{0}(x),\;\;\dfrac{\partial}{\partial
t}h_{ij}(x,0)=\varepsilon g_{ij}^{1}(x),
\end{cases}
\end{equation}
where
\begin{eqnarray*}
\bar{H}_{ij}\left(h_{kl},\dfrac{\partial h_{kl}}{\partial
x^{p}},\dfrac{\partial^{2}h_{kl}}{\partial x^{p}\partial
x^{q}}\right)&=&H^{kl}\dfrac{\partial^{2}h_{ij}}{\partial
x^{k}x^{l}}+\widetilde
H_{ij}\left(\delta_{kl}+h_{kl},\dfrac{\partial h_{kl}}{\partial
x^{p}}\right)\\\nonumber
&=&-\dfrac{1}{2}\delta^{kl}\delta^{ab}\left(\dfrac{\partial
h_{ai}}{\partial x^{k}}+\dfrac{\partial h_{ak}}{\partial
x^{i}}-\dfrac{\partial h_{ik}}{\partial
x^{a}}\right)\left(\dfrac{\partial h_{bj}}{\partial
x^{l}}+\dfrac{\partial h_{bl}}{\partial x^{j}}-\dfrac{\partial
h_{jl}}{\partial x^{b}}\right)\\\nonumber
&&-\dfrac{1}{2}\delta^{pr}\delta^{qs}\left(\dfrac{\partial
h_{ir}}{\partial x^{s}}+\dfrac{\partial h_{is}}{\partial
x^{r}}-\dfrac{\partial h_{rs}}{\partial x^{i}}\right)\dfrac{\partial
h_{pq}}{\partial
x^{j}}\\\nonumber&&-\dfrac{1}{2}\delta^{pr}\delta^{qs}\left(\dfrac{\partial
h_{jr}}{\partial x^{s}}+\dfrac{\partial h_{js}}{\partial
x^{r}}-\dfrac{\partial h_{rs}}{\partial x^{j}}\right)\dfrac{\partial
h_{pq}}{\partial x^{i}}\\\nonumber&&-h^{kl}\dfrac{\partial
^{2}h_{ij}}{\partial x^{k}\partial
x^{l}}+O\left(|h_{kl}|+\left|\dfrac{\partial h_{kl}}{\partial
x^{p}}\right|+\left|\dfrac{\partial^{2}h_{kl}}{\partial
x^{p}\partial x^{q}}\right|\right)^{3}.
\end{eqnarray*}

Let $$\hat{\lambda}=\left(h_{kl}, \dfrac{\p h_{kl}}{\p x^p},
\dfrac{\p^2 h_{kl}}{\p x^p\p x^q}\right)\quad (p,q, k,
l=1,2\cdots, n).$$ The nonlinear term
$$\bar{H}_{ij}(\hat{\lambda})=\bar{H}_{ij}\left(h_{kl},\dfrac{\partial
h_{kl}}{\partial x^{p}},\dfrac{\partial^{2}h_{kl}}{\partial
x^{p}\partial x^{q}}\right)$$ in $(4.3)$ is smooth in a
neighborhood about $\hat{\lambda}=0$ and satisfies
$$\bar{H}_{ij}(\hat{\lambda})=O\left(|\hat{\lambda}|^{2}\right)\quad (i,j=1,2,\cdots, n).$$
By the well-known global existence results for the nonlinear wave
equation (e.g., see \cite{c}, \cite{ho}, \cite{jo}, \cite{k}),
there exists a unique global smooth solution $(h_{ij}(x,t))$ for
the Cauchy problem (4.3) or (4.2). Thus, the proof of Theorem 4.1
is complete. $\quad\P$

\section{Wave character of the curvatures}
The hyperbolic geometric flow is a system of hyperbolic evolution
equations on the metrics. The evolution of the metrics implies a
system of nonlinear  wave equations for the Riemannian curvature
tensor $R_{ijkl}$, the Ricci curvature tensor $R_{ij}$ and the
scalar curvature $R$ which we will derive.

Let $\mathscr{M}$ be an $n$-dimensional complete manifold. We
consider the hyperbolic geometric flow on $\mathscr{M}$, that is,
\begin{equation}\dfrac{\p^2}{\p t^2}g_{ij}(x,t)=-2R_{ij}(x,t).
\end{equation}
We now want to find the evolution
equations for the Riemannian curvature tensor $R_{ijkl}$, the
Ricci curvature tensor $R_{ij}$ and the scalar curvature $R$.

Direct computations yield
\begin{eqnarray*}\Gamma_{jl}^h&=&\dfrac{1}{2} g^{hm}\left(\dfrac{\partial
g_{mj}}{\partial x^{l}}+\dfrac{\partial g_{ml}}{\partial
x^{j}}-\dfrac{\partial g_{jl}}{\partial x^{m}}\right),\vspace{2mm}\\
\dfrac{\p}{\p t}\Gamma_{jl}^{h}&=&\dfrac{1}{2}
g^{hm}\left(\dfrac{\partial^2 g_{mj}}{\partial x^{l}\p
t}+\dfrac{\partial^2 g_{ml}}{\partial x^{j}\p t}-\dfrac{\partial^2
g_{jl}}{\partial x^{m}\p t}\right)+ \dfrac{1}{2} \dfrac{\p
g^{hm}}{\p t}\left(\dfrac{\partial g_{mj}}{\partial
x^{l}}+\dfrac{\partial g_{ml}}{\partial x^{j}}-\dfrac{\partial
g_{jl}}{\partial x^{m}}\right),\vspace{2mm}\\
\dfrac{\p^2}{\p t^2}{\Gamma_{jl}^h}&=&\dfrac{1}{2}\dfrac{\p^2
g^{hm}}{\p t^2}\left(\dfrac{\partial g_{mj}}{\partial
x^{l}}+\dfrac{\partial g_{ml}}{\partial x^{j}}-\dfrac{\partial
g_{jl}}{\partial x^{m}}\right)+2\cdot \dfrac{1}{2}\dfrac{\p
g^{hm}}{\p t}\left(\dfrac{\partial^2 g_{mj}}{\partial x^{l}\p
t}+\dfrac{\partial^2 g_{ml}}{\partial x^{j}\p t}-\dfrac{\partial^2
g_{jl}}{\partial
x^{m}\p t}\right)\vspace{2mm}\\
&&+\dfrac{1}{2}g^{hm}\left(\dfrac{\p}{\p x^{l}}\left(\dfrac{\p^2
g_{mj}}{\p t^2}\right)+\dfrac{\p}{\p x^{j}}\left(\dfrac{\p^2
g_{ml}}{\p t^2}\right)-\dfrac{\p}{\p
x^{m}}\left(\dfrac{\p^2 g_{jl}}{\p t^2}\right)\right),\vspace{2mm}\\
R_{ijl}^{\ h}&=&\dfrac{\p \Gamma_{jl}^h}{\p
x^i}-\dfrac{\p\Gamma^{h}_{il}}{\p
x^j}+\Gamma_{ip}^h\Gamma_{jl}^p-\Gamma_{jp}^h\Gamma_{il}^p,\vspace{2mm}\\
\dfrac{\p^2}{\p t^2}R_{ijl}^{\ h}&=&\dfrac{\p }{\p
x^i}\left(\dfrac{\p^2}{\p t^2}\Gamma_{jl}^h\right)-\dfrac{\p}{\p
x^j}\left(\dfrac{\p^2}{\p
t^2}\Gamma^{h}_{il}\right)+\dfrac{\p^2}{\p
t^2}\left(\Gamma_{ip}^h\Gamma_{jl}^h-\Gamma_{jp}^h\Gamma_{il}^p\right),\vspace{2mm}\\
\dfrac{\p^2}{\p t^2}R_{ijkl}&=&\dfrac{\p^2}{\p
t^2}\left(g_{hk}R_{ijl}^{h}\right)=g_{hk}\dfrac{\p^2}{\p t^2}
R_{ijl}^{h}+R_{ijl}^{h}\dfrac{\p^2g_{hk}}{\p t^2}+2\dfrac{\p
g_{hk}}{\p t}\dfrac{\p }{\p t}R_{ijl}^h\vspace{2mm}\\
&=&g_{hk}\left[\dfrac{\p }{\p x^i}\left(\dfrac{\p^2}{\p
t^2}\Gamma_{jl}^h\right)-\dfrac{\p}{\p x^j}\left(\dfrac{\p^2}{\p
t^2}\Gamma^{h}_{il}\right)+\dfrac{\p^2}{\p
t^2}\left(\Gamma_{ip}^h\Gamma_{jl}^p-\Gamma_{jp}^h\Gamma_{il}^p\right)\right]\vspace{2mm}\\
&&+2\dfrac{\p g_{hk}}{\p t}\left[\dfrac{\p }{\p
x^i}\left(\dfrac{\p}{\p t}\Gamma_{jl}^h\right)-\dfrac{\p}{\p
x^j}\left(\dfrac{\p}{\p t}\Gamma^{h}_{il}\right)+\dfrac{\p}{\p
t}\left(\Gamma_{ip}^h\Gamma_{jl}^p-\Gamma_{jp}^h\Gamma_{il}^p\right)\right]+R_{ijl}^{h}\dfrac{\p^2g_{hk}}{\p
t^2}.\end{eqnarray*} We choose the normal coordinates around a
fixed point $p$ on $\mathscr{M}$ such that
$$\Gamma_{ij}^{k}(p)=0,$$
or, equivalently, $$\dfrac{\partial g_{ij}}{\partial
 x^k}(p)=0.$$ Then we have
\begin{equation}\begin{array}{lll}\dfrac{\p^2}{\p t^2}R_{ijkl}&=&g_{hk}\dfrac{\p}{\p
 x^i}\left[\dfrac{1}{2}\dfrac{\p^2 g^{hm}}{\p t^2}\left(\dfrac{\partial
g_{mj}}{\partial x^{l}}+\dfrac{\partial g_{ml}}{\partial
x^{j}}-\dfrac{\partial g_{jl}}{\partial x^{m}}\right)+2\cdot
\dfrac{1}{2}\dfrac{\p g^{hm}}{\p t}\left(\dfrac{\partial^2
g_{mj}}{\partial x^{l}\p t}+\dfrac{\partial^2 g_{ml}}{\partial
x^{j}\p
t}-\dfrac{\partial^2 g_{jl}}{\partial x^{m}\p t}\right)\right]\vspace{3mm}\\
&&+g_{hk}\dfrac{\p}{\p
 x^i}\left[\dfrac{1}{2}g^{hm}\left(\dfrac{\p}{\p x^{l}}\left(\dfrac{\p^2
g_{mj}}{\p t^2}\right)+\dfrac{\p}{\p x^{j}}\left(\dfrac{\p^2
g_{ml}}{\p t^2}\right)-\dfrac{\p}{\p
x^{m}}\left(\dfrac{\p^2 g_{jl}}{\p t^2}\right)\right)\right]\vspace{3mm}\\
&&-g_{hk}\dfrac{\p}{\p
 x^j}\left[\dfrac{1}{2}\dfrac{\p g^{hm}}{\p t^2}\left(\dfrac{\partial
g_{mi}}{\partial x^{l}}+\dfrac{\partial g_{ml}}{\partial
x^{i}}-\dfrac{\partial g_{il}}{\partial x^{m}}\right)+2\cdot
\dfrac{1}{2}\dfrac{\p g^{hm}}{\p t}\left(\dfrac{\partial^2
g_{mi}}{\partial x^{l}\p t}+\dfrac{\partial^2 g_{ml}}{\partial
x^{i}\p
t}-\dfrac{\partial^2 g_{il}}{\partial x^{m}\p t}\right)\right]\vspace{3mm}\\
&&-g_{hk}\dfrac{\p}{\p
 x^j}\left[\dfrac{1}{2}g^{hm}\left(\dfrac{\p}{\p x^{l}}\left(\dfrac{\p^2
g_{mi}}{\p t^2}\right)+\dfrac{\p}{\p x^{i}}\left(\dfrac{\p^2
g_{ml}}{\p t^2}\right)-\dfrac{\p}{\p
x^{m}}\left(\dfrac{\p^2 g_{il}}{\p t^2}\right)\right)\right]\vspace{3mm}\\
&&+2g_{hk}\left( \dfrac{\p}{\p t}\Gamma^{h}_{ip}\dfrac{\p}{\p
t}\Gamma^{p}_{jl}-\dfrac{\p}{\p t}\Gamma^{h}_{jp}\dfrac{\p}{\p
t}\Gamma^{p}_{il}\right)\vspace{3mm}\\
&&+2\dfrac{\p g_{hk}}{\p t}\left[\dfrac{1}{2}\dfrac{\p}{\p
x^i}\left(\dfrac{\p g^{hm}}{\p t}\left(\dfrac{\partial
g_{mj}}{\partial x^{l}}+\dfrac{\partial g_{ml}}{\partial
x^{j}}-\dfrac{\partial g_{jl}}{\partial
x^{m}}\right)\right)-\dfrac{1}{2}\dfrac{\p}{\p x^j}\left(\dfrac{\p
g^{hm}}{\p t}\left(\dfrac{\partial g_{mi}}{\partial
x^{l}}+\dfrac{\partial g_{ml}}{\partial x^{i}}-\dfrac{\partial
g_{il}}{\partial x^{m}}\right)\right)\right]\vspace{3mm}\\
&&+2\dfrac{\p g_{hk}}{\p t}\dfrac{1}{2}g^{hm}\dfrac{\p}{\p
x^i}\left(\dfrac{\p}{\p x^{l}}\left(\dfrac{\p g_{mj}}{\p
t}\right)+\dfrac{\p}{\p x^{j}}\left(\dfrac{\p g_{ml}}{\p
t}\right)-\dfrac{\p}{\p x^{m}}\left(\dfrac{\p g_{jl}}{\p t}\right)\right)\vspace{3mm}\\
&&-2\dfrac{\p g_{hk}}{\p t}\dfrac{1}{2}g^{hm}\dfrac{\p}{\p
x^j}\left(\dfrac{\p}{\p x^{l}}\left(\dfrac{\p g_{mi}}{\p
t}\right)+\dfrac{\p}{\p x^{i}}\left(\dfrac{\p g_{ml}}{\p
t}\right)-\dfrac{\p}{\p x^{m}}\left(\dfrac{\p g_{il}}{\p
t}\right)\right)+R_{ijl}^{h}\dfrac{\p^{2}g_{hk}}{\p
t^2}\vspace{3mm}\\
&=&\dfrac{1}{2}g_{hk}\dfrac{\p^2 g^{hm}}{\p t^2}\left[\dfrac{\p}{\p
x^i}\left(\dfrac{\partial g_{mj}}{\partial x^{l}}+\dfrac{\partial
g_{ml}}{\partial x^{j}}-\dfrac{\partial g_{jl}}{\partial
x^{m}}\right)-\dfrac{\p}{\p x^j}\left(\dfrac{\p g^{hm}}{\p
t}\left(\dfrac{\partial g_{mi}}{\partial x^{l}}+\dfrac{\partial
g_{ml}}{\partial x^{i}}-\dfrac{\partial g_{il}}{\partial
x^{m}}\right)\right)\right]\vspace{3mm}\\
&&+g_{hk}\dfrac{\p^2 g^{hm}}{\p x^i\p t}\left(\dfrac{\partial
g_{mj}}{\partial x^{l}\p t}+\dfrac{\partial g_{ml}}{\partial
x^{j}\p t}-\dfrac{\partial g_{jl}}{\partial x^{m}\p
t}\right)-g_{hk}\dfrac{\p^2 g^{hm}}{\p x^j\p
t}\left(\dfrac{\partial g_{mi}}{\partial x^{l}\p
t}+\dfrac{\partial g_{ml}}{\partial
x^{i}\p t}-\dfrac{\partial g_{il}}{\partial x^{m}\p t}\right)\vspace{3mm}\\
&&+g_{hk}\dfrac{\p g^{hm}}{\p t}\left[\dfrac{\p}{\p
x^i}\left(\dfrac{\partial g_{mj}}{\partial x^{l}\p
t}+\dfrac{\partial g_{ml}}{\partial x^{j}\p t}-\dfrac{\partial
g_{jl}}{\partial x^{m}\p t }\right)-\dfrac{\p}{\p
x^j}\left(\dfrac{\partial^2 g_{mi}}{\partial x^{l}\p
t}+\dfrac{\partial^2 g_{ml}}{\partial x^{i}\p t}-\dfrac{\partial^2
g_{il}}{\partial
x^{m}\p t}\right)\right]\vspace{3mm}\\
&&+\dfrac{1}{2}\left[\dfrac{\p^2}{\p x^i \p x^l}\left(\dfrac{\p^2
g_{kj}}{\p t^2}\right)+\dfrac{\p^2}{\p x^i \p
x^j}\left(\dfrac{\p^2 g_{kl}}{\p t^2}\right)-\dfrac{\p^2}{\p x^i
\p x^k}\left(\dfrac{\p^2 g_{jl}}{\p
t^2}\right)\right]\vspace{3mm}\\
&&-\dfrac{1}{2}\left[\dfrac{\p^2}{\p x^j \p x^l}\left(\dfrac{\p^2
g_{ki}}{\p t^2}\right)+\dfrac{\p^2}{\p x^j \p
x^i}\left(\dfrac{\p^2 g_{kl}}{\p t^2}\right)-\dfrac{\p^2}{\p x^j
\p x^k}\left(\dfrac{\p^2 g_{il}}{\p
t^2}\right)\right]\vspace{3mm}\\
&&+2g_{hk}\left( \dfrac{\p}{\p t}\Gamma^{h}_{ip}\dfrac{\p}{\p
t}\Gamma^{p}_{jl}-\dfrac{\p}{\p t}\Gamma^{h}_{jp}\dfrac{\p}{\p
t}\Gamma^{p}_{il}\right)\vspace{3mm}\\
&&+\dfrac{\p g_{hk}}{\p t}\dfrac{\p g^{hm}}{\p t}\left[\dfrac{\p}{\p
x^i}\left(\dfrac{\partial g_{mj}}{\partial x^{l}}+\dfrac{\partial
g_{ml}}{\partial x^{j}}-\dfrac{\partial g_{jl}}{\partial
x^{m}}\right)-\dfrac{\p}{\p x^j}\left(\dfrac{\partial
g_{mi}}{\partial x^{l}}+\dfrac{\partial g_{ml}}{\partial
x^{i}}-\dfrac{\partial
g_{il}}{\partial x^{m}}\right)\right]\vspace{3mm}\\
&&+ \dfrac{\p g_{hk}}{\p t}{ g^{hm}}\left[\dfrac{\p}{\p
x^i}\left(\dfrac{\partial g_{mj}}{\partial x^{l}\p
t}+\dfrac{\partial g_{ml}}{\partial x^{j}\p t}-\dfrac{\partial
g_{jl}}{\partial x^{m}\p t}\right)-\dfrac{\p}{\p
x^j}\left(\dfrac{\partial g_{mi}}{\partial x^{l}\p
t}+\dfrac{\partial g_{ml}}{\partial x^{i}\p t}-\dfrac{\partial
g_{il}}{\partial x^{m}\p t}\right)\right]\vspace{3mm}\\
&&+R_{ijl}^{h}\dfrac{\p^{2}g_{hk}}{\p t^2}.
\end{array}\end{equation}
Noting $g^{hm}g_{ml}=\delta^{h}_{l}$, we get
$$\begin{array}{l}\dfrac{\p g^{hm}}{\p t}=-g^{hp}g^{mq}\dfrac{\p
g_{pq}}{\p
t},\vspace{3mm}\\
\dfrac{\p^2 g^{hm}}{\p x^k\p t}=-g^{hp}g^{mq}\dfrac{\p g_{pq}}{\p
x^k\p t},\vspace{3mm}\\ \dfrac{\partial^{2} g^{hm}}{\partial t^{2}}
= -g^{hp}g^{mq}\dfrac{\partial^{2} g_{pq}}{\partial
t^{2}}+{2}g^{hp}g^{rq}g^{sm}\dfrac{\partial g_{pq}}{\partial
t}\dfrac{\partial g_{rs}}{\partial t}.
\end{array}$$
Thus, it follows from (5.2) that
\begin{equation}\begin{array}{lll}
\dfrac{\partial^{2}}{\partial
t^{2}}R_{ijkl}&=&\left(-\dfrac{1}{2}g^{pm}\dfrac{\partial^{2}
g_{kp}}{\partial t^{2}}+g^{rq}g^{pm}\dfrac{\partial g_{kq}}{\partial
t}\dfrac{\partial g_{rp}}{\partial
t}\right)\vspace{2mm}\\&&\times\left[\dfrac{\partial}{\partial
x^{i}}\left(\dfrac{\partial g_{mj}}{\partial x^{l}}+\dfrac{\partial
g_{ml}}{\partial x^{j}}-\dfrac{\partial g_{jl}}{\partial
x^{m}}\right)- \dfrac{\partial}{\partial x^{j}}\left(\dfrac{\partial
g_{mi}}{\partial x^{l}}+\dfrac{\partial g_{ml}}{\partial
x^{i}}-\dfrac{\partial g_{il}}{\partial
x^{m}}\right)\right]\vspace{2mm}\\&&-g^{pm}\dfrac{\partial^{2}
g_{kp}}{\partial x^{i}\partial t}\left(\dfrac{\partial^{2}
g_{mj}}{\partial x^{l}\partial t}+\dfrac{\partial^{2}
g_{ml}}{\partial x^{j}\partial t}-\dfrac{\partial^{2}
g_{jl}}{\partial x^{m}\partial t}\right) + g^{pm}\dfrac{\partial^{2}
g_{kp}}{\partial x^{j}\partial t}\left(\dfrac{\partial^{2}
g_{mi}}{\partial x^{l}\partial t}+\dfrac{\partial^{2}
g_{ml}}{\partial x^{i}\partial t}-\dfrac{\partial^{2}
g_{il}}{\partial x^{m}\partial
t}\right)\vspace{2mm}\\&&+\dfrac{1}{2}\left[\dfrac{\partial^{2}}{\partial
x^{i}\partial x^{l}}\left(\dfrac{\partial^{2} g_{kj}}{\partial
t^{2}}\right)+\dfrac{\partial^{2}}{\partial x^{i}\partial
x^{j}}\left(\dfrac{\partial^{2} g_{kl}}{\partial t^{2}}\right)-
\dfrac{\partial^{2}}{\partial x^{i}\partial
x^{k}}\left(\dfrac{\partial^{2} g_{jl}}{\partial
t^{2}}\right)\right]\vspace{2mm}\\&&-\dfrac{1}{2}\left[\dfrac{\partial^{2}}{\partial
x^{j}\partial x^{l}}\left(\dfrac{\partial^{2} g_{ki}}{\partial
t^{2}}\right)+\dfrac{\partial^{2}}{\partial x^{j}\partial
x^{i}}\left(\dfrac{\partial^{2} g_{kl}}{\partial
t^{2}}\right)-\dfrac{\partial^{2}}{\partial x^{j}\partial
x^{k}}\left(\dfrac{\partial^{2} g_{il}}{\partial
t^{2}}\right)\right]\vspace{2mm}\\&&+
{2}g_{hk}\left(\dfrac{\partial}{\partial
t}\Gamma^{h}_{ip}\cdot\dfrac{\partial}{\partial
t}\Gamma^{p}_{jl}-\dfrac{\partial}{\partial
t}\Gamma^{h}_{jp}\cdot\dfrac{\partial}{\partial
t}\Gamma^{p}_{il}\right)+R^h_{ijl}\dfrac{\partial^{2}g_{hk}}{\partial
t^{2}}\vspace{2mm}\\&&-g^{hp}g^{mq}\dfrac{\partial g_{nk}}{\partial
t}\dfrac{\partial g_{pq}}{\partial t}\left[\dfrac{\partial}{\partial
x^{i}}\left(\dfrac{\partial g_{mj}}{\partial x^{l}}+\dfrac{\partial
g_{ml}}{\partial x^{j}}-\dfrac{\partial g_{jl}}{\partial
x^{m}}\right)- \dfrac{\partial}{\partial x^{j}}\left(\dfrac{\partial
g_{mi}}{\partial x^{l}}+\dfrac{\partial g_{ml}}{\partial
x^{i}}-\dfrac{\partial g_{il}}{\partial
x^{m}}\right)\right]\vspace{2mm}\\&=& -\dfrac{\partial^{2}
g_{kp}}{\partial t^{2}}\left(\dfrac{\partial}{\partial
x^{i}}\Gamma^{p}_{jl}-\dfrac{\partial}{\partial
x^{j}}\Gamma^{p}_{il}\right)+ {2}g^{rq}\dfrac{\partial
g_{kq}}{\partial t} \dfrac{\partial g_{rp}}{\partial
t}\left(\dfrac{\partial}{\partial
x^{i}}\Gamma^{p}_{jl}-\dfrac{\partial}{\partial
x^{j}}\Gamma^{p}_{il}\right)\vspace{2mm}\\&&- \dfrac{\partial^{2}
g_{kp}}{\partial x^{i}\partial t}\cdot
g^{pm}\left(\dfrac{\partial^{2} g_{mj}}{\partial x^{l}\partial
t}+\dfrac{\partial^{2} g_{ml}}{\partial x^{j}\partial
t}-\dfrac{\partial^{2} g_{jl}}{\partial x^{m}\partial
t}\right)+\dfrac{\partial^{2} g_{kp}}{\partial x^{j}\partial t}\cdot
g^{pm}\left(\dfrac{\partial^{2} g_{mi}}{\partial x^{l}\partial
t}+\dfrac{\partial^{2} g_{ml}}{\partial x^{i}\partial
t}-\dfrac{\partial^{2} g_{il}}{\partial x^{m}\partial
t}\right)\vspace{2mm}\\&&+\dfrac{1}{2}\left[\dfrac{\partial^{2}}{\partial
x^{i}\partial
x^{l}}\left(-{2}R_{kj}\right)+\dfrac{\partial^{2}}{\partial
x^{i}\partial x^{j}}\left(-{2}R_{kl}\right)-
\dfrac{\partial^{2}}{\partial x^{i}\partial
x^{k}}\left(-{2}R_{jl}\right)\right]\vspace{2mm}\\&&-\dfrac{1}{2}\left[\dfrac{\partial^{2}}{\partial
x^{j}\partial
x^{l}}\left(-{2}R_{ik}\right)+\dfrac{\partial^{2}}{\partial
x^{i}\partial
x^{j}}\left(-{2}R_{kl}\right)-\dfrac{\partial^{2}}{\partial
x^{j}\partial x^{k}}\left(-{2}R_{il}\right)\right]\vspace{2mm}\\&&+
{2}g_{hk}\left[\dfrac{\partial}{\partial
t}\Gamma^{h}_{ip}\cdot\dfrac{\partial}{\partial
t}\Gamma^{p}_{jl}-\dfrac{\partial}{\partial
t}\Gamma^{h}_{jp}\cdot\dfrac{\partial}{\partial
t}\Gamma^{p}_{il}\right]-\dfrac{\partial g_{hk}}{\partial
t}\dfrac{\partial g_{pq}}{\partial
t}g^{hp}\cdot{2}\left(\dfrac{\p}{\p x_i}\Gamma^q_{jl}-\dfrac{\p}{\p
x_j}\Gamma^q_{il}
\right)+R^h_{ijl}\dfrac{\partial^{2}g_{hk}}{\partial
t^{2}}\vspace{2mm}\\&=&
\dfrac{1}{2}\left[\dfrac{\partial^{2}}{\partial x^{i}\partial
x^{l}}\left(-{2}R_{kj}\right)+\dfrac{\partial^{2}}{\partial
x^{i}\partial
x^{j}}\left(-{2}R_{kl}\right)-\dfrac{\partial^{2}}{\partial
x^{i}\partial
x^{k}}\left(-{2}R_{jl}\right)\right]\vspace{2mm}\\&&-\dfrac{1}{2}\left[\dfrac{\partial^{2}}{\partial
x^{j}\partial x^{l}}\left(-{2}R_{ki}\right)+
\dfrac{\partial^{2}}{\partial x^{i}\partial
x^{j}}\left(-{2}R_{kl}\right)-\dfrac{\partial^{2}}{\partial
x^{j}\partial
x^{k}}\left(-{2}R_{il}\right)\right]\vspace{2mm}\\&&-g^{pm}\dfrac{\partial^{2}
g_{kp}}{\partial x^{i}\partial t}\left(\dfrac{\partial^{2}
g_{mj}}{\partial x^{l}\partial t}+\dfrac{\partial^{2}
g_{ml}}{\partial x^{j}\partial t}-\dfrac{\partial^{2}
g_{jl}}{\partial x^{m}\partial t}\right)+ g^{pm}\dfrac{\partial^{2}
g_{kp}}{\partial x^{j}\partial t}\left(\dfrac{\partial^{2}
g_{mi}}{\partial x^{l}\partial t}+\dfrac{\partial^{2}
g_{ml}}{\partial x^{i}\partial t}-\dfrac{\partial^{2}
g_{il}}{\partial x^{m}\partial
t}\right)\vspace{2mm}\\&&+{2}g_{hk}\left(\dfrac{\partial}{\partial
t}\Gamma^{h}_{ip}\cdot\dfrac{\partial}{\partial
t}\Gamma^{p}_{jl}-\dfrac{\partial}{\partial
t}\Gamma^{h}_{jp}\cdot\dfrac{\partial}{\partial
t}\Gamma^{p}_{il}\right).
\end{array}\end{equation}
On the one hand, we have
$$\dfrac{\partial^{2}}{\partial x^{i} \partial
x^{l}}R_{jk}=\nabla_{i}\nabla_{l}R_{jk}+\nabla_{i}\Gamma^{p}_{lk}\cdot
R_{jp}+\nabla_{i}\Gamma^{p}_{lj}\cdot R_{kp}.$$ Then
\begin{eqnarray}\begin{array}{lll}
&&\dfrac{1}{2}\left[\dfrac{\partial^{2}}{\partial x^{i}\partial
x^{l}}\left(-{2}R_{kj}\right)+\dfrac{\partial^{2}}{\partial
x^{i}\partial
x^{j}}\left(-{2}R_{kl}\right)-\dfrac{\partial^{2}}{\partial
x^{i}\partial
x^{k}}\left(-{2}R_{jl}\right)\right]\vspace{2mm}\\&&-\dfrac{1}{2}\left[\dfrac{\partial^{2}}{\partial
x^{j}\partial x^{l}}\left(-{2}R_{ki}\right)+
\dfrac{\partial^{2}}{\partial x^{i}\partial
x^{j}}\left(-{2}R_{kl}\right)-\dfrac{\partial^{2}}{\partial
x^{j}\partial x^{k}}\left(-{2}R_{il}\right)\right]\vspace{2mm}\\&=&
-\nabla_{i}\nabla_{l} R_{kj}-\nabla_{i}\Gamma^{p}_{lk}
R_{jp}-\nabla_{i}\Gamma^{p}_{lj} R_{kp}- \nabla_{i}\nabla_{j}
R_{kl}-\nabla_{i}\Gamma^{p}_{jk} R_{lp}-\nabla_{i}\Gamma^{p}_{jl}
R_{kp}\vspace{2mm}\\&&+\nabla_{i}\nabla_{k}
R_{jl}+\nabla_{i}\Gamma^{p}_{kj} R_{lp}+\nabla_{i}\Gamma^{p}_{kl}
R_{jp}+\nabla_{i}\nabla_{l} R_{ki}+ \nabla_{j}\Gamma^{p}_{kl}
R_{ip}+\nabla_{j}\Gamma^{p}_{li}
R_{kp}\vspace{2mm}\\&&+\nabla_{j}\nabla_{i}
R_{kl}+\nabla_{j}\Gamma^{p}_{ik} R_{lp}+\nabla_{j}\Gamma^{p}_{il}
R_{kp}-\nabla_{j}\nabla_{k} R_{il}-\nabla_{j}\Gamma^{p}_{ki} R_{pl}-
\nabla_{j}\Gamma^{p}_{kl} R_{pi}\vspace{2mm}\\&=&
-\nabla_{i}\nabla_{l} R_{kj}-\nabla_{i}\nabla_{j}
R_{kl}+\nabla_{i}\nabla_{k} R_{jl}+\nabla_{j}\nabla_{l}
R_{ki}+\nabla_{j}\nabla_{i} R_{ki}-\nabla_{j}\nabla_{k}
R_{il}\vspace{2mm}\\&&+
R_{ip}\left(\nabla_{j}\Gamma^{p}_{lk}-\nabla_{j}\Gamma^{p}_{kl}\right)+
R_{jp}\left(-\nabla_{i}\Gamma^{p}_{lk}+\nabla_{i}\Gamma^{p}_{kl}\right)\vspace{2mm}\\&&+
R_{kp}\left(-\nabla_{i}\Gamma^{p}_{lj}-\nabla_{i}\Gamma^{p}_{jl}+
\nabla_{j}\Gamma^{p}_{li}+\nabla_{j}\Gamma^{p}_{il}\right)+
R_{lp}\left(-\nabla_{i}\Gamma^{p}_{jk}+\nabla_{i}\Gamma^{p}_{kj}+
\nabla_{j}\Gamma^{p}_{ik}-\nabla_{j}\Gamma^{p}_{ki}\right)\vspace{2mm}\\&=&
-\nabla_{i}\nabla_{l}R_{kj}+\nabla_{i}\nabla_{k}R_{jl}+\nabla_{j}\nabla_{l}R_{ki}-
\nabla_{j}\nabla_{k}R_{il}\vspace{2mm}\\
&&-R_{ijlp}g^{pq}R_{qk}-R_{ijkp}g^{pq}R_{ql}+R_{kp}\left(-{2}R_{ijl}^p\right)\vspace{2mm}\\&=&
-\nabla_{i}\nabla_{l}R_{kj}+
\nabla_{i}\nabla_{k}R_{jl}+\nabla_{j}\nabla_{l}R_{ki}-\nabla_{j}\nabla_{k}R_{il}-
g^{pq}\left(R_{ijlp}R_{kq}+R_{ijkp}R_{lq}\right).\end{array}\end{eqnarray}
On the other hand, we have
\begin{eqnarray}\begin{array}{lll}
&&-g^{pm}\dfrac{\partial^{2} g_{kp}}{\partial x^{i}\partial
t}\left(\dfrac{\partial^{2} g_{mj}}{\partial x^{l}\partial
t}+\dfrac{\partial^{2} g_{ml}}{\partial x^{j}\partial
t}-\dfrac{\partial^{2} g_{jl}}{\partial x^{m}\partial t}\right)+
g^{pm}\dfrac{\partial^{2} g_{kp}}{\partial x^{j}\partial
t}\left(\dfrac{\partial^{2} g_{mi}}{\partial x^{l}\partial
t}+\dfrac{\partial^{2} g_{ml}}{\partial x^{i}\partial
t}-\dfrac{\partial^{2} g_{il}}{\partial x^{m}\partial
t}\right)\vspace{2mm}\\&&+ {2}g_{hk}\left(\dfrac{\partial}{\partial
t}\Gamma^{h}_{ip}\cdot\dfrac{\partial}{\partial
t}\Gamma^{p}_{jl}-\dfrac{\partial}{\partial
t}\Gamma^{h}_{jp}\cdot\dfrac{\partial}{\partial
t}\Gamma^{p}_{il}\right)\vspace{2mm}\\&=&-g^{pm}\dfrac{\partial^{2}
g_{kp}}{\partial x^{i}\partial t}\left(\dfrac{\partial^{2}
g_{mj}}{\partial x^{l}\partial t}+\dfrac{\partial^{2}
g_{ml}}{\partial x^{j}\partial t}-\dfrac{\partial^{2}
g_{jl}}{\partial x^{m}\partial t}\right)+ g^{pm}\dfrac{\partial^{2}
g_{kp}}{\partial x^{j}\partial t}\left(\dfrac{\partial^{2}
g_{mi}}{\partial x^{l}\partial t}+\dfrac{\partial^{2}
g_{ml}}{\partial x^{i}\partial t}-\dfrac{\partial^{2}
g_{il}}{\partial x^{m}\partial
t}\right)\vspace{2mm}\\&&+\dfrac{1}{2}g^{pm}\left(\dfrac{\partial^{2}
g_{ki}}{\partial x^{p}\partial t}+\dfrac{\partial^{2}
g_{kp}}{\partial x^{i}\partial t}-\dfrac{\partial^{2}
g_{ip}}{\partial x^{k}\partial t}\right)\left(\dfrac{\partial^{2}
g_{mj}}{\partial x^{l}\partial t}+ \dfrac{\partial^{2}
g_{ml}}{\partial x^{j}\partial t}-\dfrac{\partial^{2}
g_{jl}}{\partial x^{m}\partial
t}\right)\vspace{2mm}\\&&-\dfrac{1}{2}g^{pm}\left(\dfrac{\partial^{2}
g_{kj}}{\partial x^{p}\partial t}+\dfrac{\partial^{2}
g_{kp}}{\partial x^{j}\partial t}-\dfrac{\partial^{2}
g_{jp}}{\partial x^{k}\partial t}\right) \left(\dfrac{\partial^{2}
g_{mi}}{\partial x^{l}\partial t}+\dfrac{\partial^{2}
g_{ml}}{\partial x^{i}\partial t}-\dfrac{\partial^{2}
g_{il}}{\partial x^{m}\partial t}\right)\vspace{2mm}\\&=&
 g^{pm}\left\{{\dfrac{1}{2}}\left(\dfrac{\partial^{2} g_{mj}}{\partial
x^{l}\partial t}+\dfrac{\partial^{2} g_{ml}}{\partial
x^{j}\partial
 t}-\dfrac{\partial^{2} g_{jl}}{\partial x^{m}\partial
 t}\right)\left(\dfrac{\partial^{2} g_{ki}}{\partial x^{p}\partial
 t}-\dfrac{\partial^{2} g_{ip}}{\partial x^{k}\partial
 t}-\dfrac{\partial^{2} g_{kp}}{\partial x^{i}\partial
 t}\right)\right.\vspace{2mm}\\&&\left.-\dfrac{1}{2}\left(\dfrac{\partial^{2} g_{mi}}{\partial x^{l}\partial
 t}+\dfrac{\partial^{2} g_{ml}}{\partial x^{i}\partial
 t}-\dfrac{\partial^{2} g_{il}}{\partial x^{m}\partial
 t}\right) \left(\dfrac{\partial^{2} g_{kj}}{\partial x^{p}\partial t}-
 \dfrac{\partial^{2} g_{jp}}{\partial x^{k}\partial
 t}-\dfrac{\partial^{2} g_{kp}}{\partial x^{j}\partial t}\right)\right\}\vspace{2mm}\\&=&
 2g_{pq}\left\{\dfrac{1}{2}g^{pr}\dfrac{\partial}{\partial t}\left(\dfrac{\partial g_{rj}}
 {\partial x^l}+\dfrac{\partial g_{rl}}{\partial x^j}
 -\dfrac{\partial g_{jl}}{\partial x^r}\right)\cdot\dfrac{1}{2}g^{qs}\dfrac{\partial}{\partial t}
 \left(\dfrac{\partial g_{ki}}{\partial x^s}-\dfrac{\partial g_{is}}{\partial x^k}
 -\dfrac{\partial g_{ks}}{\partial x^i}\right)\right.\vspace{2mm}\\&&\left.-\dfrac{1}{2}g^{pr}\dfrac{\partial}
 {\partial t}\left(\dfrac{\partial g_{ri}}{\partial x^l}+\dfrac{\partial g_{rl}}{\partial x^i}
 -\dfrac{\partial g_{il}}{\partial x^r}\right)\cdot\dfrac{1}{2}g^{qs}\dfrac{\partial}{\partial t}
 \left(\dfrac{\partial g_{kj}}{\partial x^s}-\dfrac{\partial g_{js}}{\partial x^k}
 -\dfrac{\partial g_{ks}}{\partial x^j}\right)\right\}\vspace{2mm}\\&=&2g_{pq}\left(\dfrac{\partial}{\partial
t}\Gamma^{p}_{il}\cdot\dfrac{\partial}{\partial
t}\Gamma^{q}_{jk}-\dfrac{\partial}{\partial
t}\Gamma^{p}_{jl}\cdot\dfrac{\partial}{\partial
t}\Gamma^{q}_{ik}\right).\end{array}\end{eqnarray} Therefore, it
follows from (5.3), (5.4) and (5.5) that
\begin{align}
\dfrac{\partial^2}{\partial t^2}R_{ijkl}= &
-\nabla_{i}\nabla_{l}R_{kj}+\nabla_{i}\nabla_{k}R_{jl}+
\nabla_{j}\nabla_{l}R_{ki}-\nabla_{j}\nabla_{k}R_{il}
-g^{pq}\left(R_{ijql}R_{kp}+R_{ijkq}R_{kp}\right)\nonumber\\
& +2g_{pq}\left(\dfrac{\partial}{\partial
t}\Gamma^{p}_{il}\cdot\dfrac{\partial}{\partial
t}\Gamma^{q}_{jk}-\dfrac{\partial}{\partial
t}\Gamma^{p}_{jl}\cdot\dfrac{\partial}{\partial
t}\Gamma^{q}_{ik}\right)\nonumber.
\end{align}
Similar to Hamilton \cite{h}, we have
\begin{Theorem} Under the hyperbolic geometric flow (5.1), the Riemannian
curvature tensor $R_{ijkl}$ satisfies the evolution equation
\begin{eqnarray}\begin{array}{lll}
\dfrac{\partial^2}{\partial t^2}R_{ijkl}= & \triangle
R_{ijkl}+2\left(B_{ijkl}-B_{ijlk}-B_{iljk}+B_{ikjl}\right)\\ &
-g^{pq}\left(R_{pjkl}R_{qi}+R_{ipkl}R_{qj}+R_{ijpl}R_{qk}+R_{ijkp}R_{ql}\right)\\
& +2g_{pq}\left(\dfrac{\partial}{\partial
t}\Gamma^{p}_{il}\cdot\dfrac{\partial}{\partial
t}\Gamma^{q}_{jk}-\dfrac{\partial}{\partial
t}\Gamma^{p}_{jl}\cdot\dfrac{\partial}{\partial
t}\Gamma^{q}_{ik}\right),\end{array}
\end{eqnarray}
where $B_{ijkl}=g^{pr}g^{qs}R_{piqj}R_{rksl}$ and $\triangle$ is
the Laplacian with respect to the evolving metric.
\end{Theorem}

\begin{Remark}

In Theorem 5.1 and Theorem 5.2 below, the term
$2g_{pq}\left(\dfrac{\partial}{\partial
t}\Gamma^{p}_{il}\cdot\dfrac{\partial}{\partial
t}\Gamma^{q}_{jk}-\dfrac{\partial}{\partial
t}\Gamma^{p}_{jl}\cdot\dfrac{\partial}{\partial
t}\Gamma^{q}_{ik}\right)$ can be written in the covariant form.
For the sake of simplicity, we omit it.
\end{Remark}

For the Ricci curvature tensor, we have
\begin{eqnarray*}
\dfrac{\partial^2}{\partial t^2}R_{ik} & =&
\dfrac{\partial^2}{\partial
t^2}\left(R_{ijkl}  g^{jl}\right)\\
& =& g^{jl}\dfrac{\partial^2}{\partial
t^2}R_{ijkl}+2\dfrac{\partial}{\partial
t}g^{jl}\cdot\dfrac{\partial}{\partial
t}R_{ijkl}+R_{ijkl}\dfrac{{\partial^2}{g^{jl}}}{\partial t^2}\\
& =&g^{jl}\dfrac{\partial^2}{\partial
t^2}R_{ijkl}-2g^{jp}g^{lq}\dfrac{\partial g_{pq}}{\partial
t}\dfrac{\partial}{\partial
t}R_{ijkl}-g^{jp}g^{lq}\dfrac{{\partial^2}g_{pq}}{\partial
t^2}R_{ijkl}+2g^{jp}g^{rq}g^{sl}\dfrac{\partial g_{pq}}{\partial
t}\dfrac{\partial g_{rs}}{\partial t} R_{ijkl}.
\end{eqnarray*}
Thus, we obtain
\begin{Theorem} Under the hyperbolic geometric flow (5.1), the
Ricci curvature tensor satisfies
\begin{equation}\begin{array}{lll}
\dfrac{\partial^2}{\partial t^2}R_{ik} & = & \triangle
R_{ik}+2g^{pr}g^{qs}R_{piqk}R_{rs}-2g^{pq}R_{pi}R_{qk}\vspace{2mm}\\
& & +2g^{jl}g_{pq}\left(\dfrac{\partial}{\partial
t}\Gamma^{p}_{il}\dfrac{\partial}{\partial
t}\Gamma^{q}_{jk}-\dfrac{\partial}{\partial
t}\Gamma^{p}_{jl}\dfrac{\partial}{\partial
t}\Gamma^{q}_{ik}\right)\vspace{2mm}\\
& & -2g^{jp}g^{lq}\dfrac{\partial g_{pq}}{\partial
t}\dfrac{\partial}{\partial
t}R_{ijkl}+2g^{jp}g^{rq}g^{sl}\dfrac{\partial g_{pq}}{\partial
t}\dfrac{\partial g_{rs}}{\partial
t}R_{ijkl}.\end{array}\end{equation}
\end{Theorem}

For the scalar curvature, we have
\begin{eqnarray*}
\dfrac{\partial^2}{\partial t^2}R & =&\dfrac{\partial^2}{\partial
t^2}\left(g^{ik} R_{ik}\right)\vspace{2mm}\\
& =&g^{ik}\dfrac{\partial^2}{\partial t^2}R_{ik}
+2\dfrac{\partial}{\partial t}R_{ik}\cdot\dfrac{\partial}{\partial
t}g^{ik}+R_{ik}\dfrac{\partial^2 g^{ik}}{\partial t^2}\vspace{2mm}\\
& =&g^{ik}\dfrac{\partial^2}{\partial
t^2}R_{ik}-2g^{ip}g^{kq}\dfrac{\partial{g_{pq}}}{\partial
t}\dfrac{\partial {R_{ik}}}{\partial
t}+R_{ik}\left(-g^{ip}g^{kq}\dfrac{\partial^2{g_{pq}}}{\partial
t^2}+2g^{ip}g^{rq}g^{sk}\dfrac{\partial g_{pq}}{\partial
t}\dfrac{\partial g_{rs}}{\partial t}\right).
\end{eqnarray*}
On the other hand,
\begin{eqnarray*}
&& 2g^{ik}g^{jl}g_{pq}\left(\dfrac{\partial}{\partial
t}\Gamma^{p}_{il}\dfrac{\partial}{\partial
t}\Gamma^{q}_{jk}-\dfrac{\partial}{\partial
t}\Gamma^{p}_{jl}\dfrac{\partial}{\partial
t}\Gamma^{q}_{ik}\right)\vspace{2mm}\\
&=&\frac{3}{2}g^{ik}g^{jl}g^{rs}\nabla_r(\dfrac{\p g_{ij}}{\p
t})\nabla_s(\dfrac{\p g_{kl}}{\p
t})-\frac{1}{2}g^{rs}\nabla_r(g^{ik}\dfrac{\p g_{ik}}{\p
t})\nabla_s(g^{jl}\dfrac{\p g_{jl}}{\p t})\vspace{2mm}\\
& &+2g^{rs}g^{jl}\nabla_r(g^{ik}\dfrac{\p g_{ik}}{\p
t})\nabla_l(\dfrac{\p g_{js}}{\p
t})-g^{ik}g^{jl}g^{rs}\nabla_r(\dfrac{\p g_{ij}}{\p
t})\nabla_l(\dfrac{\p g_{ks}}{\p
t})-2g^{ik}g^{jl}g^{rs}\nabla_i(\dfrac{\p g_{kr}}{\p
t})\nabla_j(\dfrac{\p g_{ls}}{\p t}).
\end{eqnarray*}
Then, we get
\begin{Theorem} Under the hyperbolic geometric flow (5.1), the
scalar curvature satisfies
\begin{eqnarray}\begin{array}{lll}
\dfrac{\partial^2}{\partial t^2}R&= & \triangle R+2|{\rm Ric}|^2\vspace{2mm}\\
&& +\frac{3}{2}g^{ik}g^{jl}g^{rs}\nabla_r(\dfrac{\p g_{ij}}{\p
t})\nabla_s(\dfrac{\p g_{kl}}{\p
t})-\frac{1}{2}g^{rs}\nabla_r(g^{ik}\dfrac{\p g_{ik}}{\p
t})\nabla_s(g^{jl}\dfrac{\p g_{jl}}{\p t})\vspace{2mm}\\
& &+2g^{rs}g^{jl}\nabla_r(g^{ik}\dfrac{\p g_{ik}}{\p
t})\nabla_l(\dfrac{\p g_{js}}{\p
t})-g^{ik}g^{jl}g^{rs}\nabla_r(\dfrac{\p g_{ij}}{\p
t})\nabla_l(\dfrac{\p g_{ks}}{\p
t})\vspace{2mm}\\
& &-2g^{ik}g^{jl}g^{rs}\nabla_i(\dfrac{\p g_{kr}}{\p
t})\nabla_j(\dfrac{\p g_{ls}}{\p
t})-2g^{ik}g^{jp}g^{lq}\dfrac{\partial g_{pq}}{\partial
t}\dfrac{\partial}{\partial t}R_{ijkl}\vspace{2mm}\\
& &-2g^{ip}g^{kq}\dfrac{\partial g_{pq}}{\partial t}\dfrac{\partial
R_{ik}}{\partial t}+4R_{ik}g^{ip}g^{rq}g^{sk}\dfrac{\partial
g_{pq}}{\partial t}\dfrac{\partial g_{rs}}{\partial t}.
\end{array}\end{eqnarray}
\end{Theorem}

Theorems 5.1-5.3 show that the curvatures of the hyperbolic
geometric flow possess the wave character. We will apply techniques
from hyperbolic equations to the above wave equations of curvatures
to derive various geometric results.

\section{Discussions}

The hyperbolic geometric flow describes the wave character of the
metrics and curvatures of manifolds. Many hyperbolic systems in
nature provide natural singular sets, the typical example is the
Einstein equations in general relativity which form a hyperbolic
system with a well-posed Cauchy problem. If one starts with smooth
initial data, one may end up with a singular space-time. One of
the most challenging problems is to describe the kind of natural
singularity. The famous {\it cosmic censorship} conjecture due to
Penrose is an attempt to describe such singularities (see
\cite{p}). In Kong and Liu \cite{kl}, we construct some exact
solutions of the hyperbolic geometric flow, these solutions
possess the singularities which are nothing but those described by
Penrose's conjecture. By these examples, we believe that the
hyperbolic geometric flow is a very natural and powerful tool to
understand the singularities in the nature, in particular, the
singularity described by Penrose cosmic censorship conjecture.

The Einstein equations play an essential role in general relativity.
Consider a space-time with Lorentzian metric
\begin{equation}
ds^2=g_{\mu\nu}dx^{\mu}dx^{\nu}\quad (\mu,\nu=0,1,2\cdots,n).
\end{equation}
The vacuum Einstein equations read
\begin{equation}
G_{\mu\nu}=0,
\end{equation}
where $G_{\mu\nu}$ is the Einstein tensor. We now consider the
following metric with orthogonal time-axis
\begin{equation}
ds^{2}=-dt^{2}+g_{ij}(x,t)dx^{i}dx^{j}.
\end{equation}
Substituting (6.3) into (6.2), we can obtain the equations
satisfied by the metric $g_{ij}$
\begin{equation}
\dfrac{\partial^{2}g_{ij}}{\partial
               t^{2}}=-2R_{ij}-\dfrac{1}{2}g^{pq}\dfrac{\partial g_{ij}}{\partial
               t}\dfrac{\partial g_{pq}}{\partial t}+g^{pq}\dfrac{\partial g_{ip}}{\partial
               t}\dfrac{\partial g_{jq}}{\partial t}.
\end{equation}
Neglecting the lower order terms gives the hyperbolic geometric flow
(1.1). Therefore, in this sense, the hyperbolic geometric flow can
be viewed as the leading terms in the vacuum Einstein equations with
respect to the metric (6.3). Since the hyperbolic geometric flow
only contains the main terms in the Einstein equations, it not only
becomes simpler and more symmetric, but also possesses rich and
beautiful geometric properties. In particular, in mathematics, its
Cauchy problem is well-posed and easier to handle some fundamental
problems such as the global existence and formation of
singularities; on the other hand, it can be applied to re-understand
the singularity of the universe and other important problems in
physics and cosmology (see \cite{shu}). We also believe that there
should be some relations between the solutions of the Einstein
equations and the corresponding hyperbolic geometric flows. On the
other hand, from the above discussions we have seen that the
hyperbolic geometric flow also possesses many beautiful features
similar to those of the Ricci flow, and some of the techniques in
the study of the Ricci flow can be directly used to understand the
hyperbolic geometric flow. The deep study on the hyperbolic
geometric flow may open a new way to understand the complicated
Einstein equations.

It is well known, in general relativity there is a constraint system
of equations involving an asymptotically flat metric tensor and
another symmetric tensor. There are four constraint equations and it
is therefore over-determined. Unlike this, since the time axis is
orthogonal to other space axes, the hyperbolic geometric flow does
not need to satisfy any additional constraint. More precisely, for
the Cauchy problem of the hyperbolic geometric flow, in order to
determine the solution we need two initial conditions: one is the
metric flow itself $g_{ij}(x,0)$, another is its derivative
$\frac{\partial g_{ij}}{\partial t}(x,0)$, since the time axis is
orthogonal to other axes, these initial data do not need to satisfy
any additional constraint, and therefore it is a determined system.
This is another main new feature of the hyperbolic geometric flow.

Many mathematicians, for example Shatah et al \cite{s1}-\cite{s3},
have investigated the Cauchy problem for some geometric wave
equations. The model at hand is the harmonic map problem, which is
the study of maps from the Minkowski space-time into complete
Riemannian manifolds. This kind of geometric wave equations is a
system of partial differential equations of second order, which is
the Euler-Lagrange equations of the action integral of the
harmonic map. It satisfies certain linear matching condition, and
then under suitable assumptions, has a unique small smooth
solution for all time, and possesses some interesting (decay,
energy and regularity) estimates. On the other hand, the
hyperbolic geometric flow is determined by the Ricci curvatures of
a family of Riemannian metrics on the manifold under
consideration. That is to say, the hyperbolic geometric flow
possesses itself intrinsic geometric structure and can be used to
describe the wave character of metrics and curvatures. This is
essentially different from the above harmonic map problem.

As well-known, one can understand the heat kernel from the kernel
of wave equation. This indicates that we should be able to derive
various information of the Ricci flow from that of the hyperbolic
geometric flow. Therefore it is also interesting to understand the
relations between the hyperbolic geometric flow and the Ricci
flow, the singularities of its solutions and its relation with the
geometrization theorem. This will be another interesting topic in
the sequel.

\vskip 5mm

\noindent{\Large {\bf Acknowledgements.}}
%Some parts of the work
%was completed while one of the authors, D.-X. Kong, was visiting
%the Courant Institute and Harvard University. Kong thanks
%Professors C. S. Morawetz and S.-T. Yau for their kind invitations
%and warm encouragements.
The work of Dai was supported by the fund
of post-doctor of China (Grant No. 20060401063); the work of Kong
was supported in part by the NNSF of China (Grant No. 10671124)
and the NCET of China (Grant No. NCET-05-0390); the work of Liu
was supported in part by the NSF and NSF of China.


\begin{thebibliography}{99}

\bibitem{cz} Huaidong Cao and Xiping Zhu, {\it A complete proof of the
Poincar$\acute{e}$ and geometrization conjectures - application of
the Hamilton-Perelman theory of the Ricci flow}, Asian J. Math.
{\bf 10} (2006), 165-492.

%\bibitem{ckl1} Jiazhong Chen, De-Xing Kong and Kefeng Liu,
%{\it New geometric notion redefines the big bang},
%to appear.

%\bibitem{ckl2} Jiazhong Chen, De-Xing Kong and Kefeng Liu,
%{\it Applications of hyperbolic geometric flow to cosomolgy}, in
%preparation.

\bibitem{c} D. Christodoulou, {\it Global solutions of nonlinear hyperbolic equations
for small initial data}, Comm. Pure Appl. Math. {\bf 39} (1986),
367-282.

\bibitem{fm} A.E. Fischer and J.E. Marsden, {\it The Einstein evolution equations as
a first-order quasi-linear symmetric system hyperbolic system I},
Commun. Math. Phys. {\bf 28} (1972), 1-38.

\bibitem{fo} V. Fock, The theory of space, time and gravitation, second revised
edition. Translated from the Russian by N. Kemmer. A Pergamon Press
Book The Macmillan Co., New York, 1964.

\bibitem{f} K.O. Friedrich, {\it Symmetric hyperbolic linear differential equations},
Comm. Pure Appl. Math. {\bf 7} (1954), 345-392.

\bibitem{h}  R. Hamilton, {\it Three-manifolds with positive Ricci curvature},
J. Differential Geom. {\bf 17} (1982), 255-306.

\bibitem{ho} L. H$\ddot{o}$rmander, Lectures on nonlinear hyperbolic differential
equations, Math$\acute{e}$matiques \& Applications {\bf 26},
Springer-Verlag, Berlin, 1997.

\bibitem{jo} F. John, {\it Delayed singularity formation in solutions of nonlinear
wave equations in higher dimensions}, Comm. Pure Appl. Math. {\bf
29} (1976), 649-682.

\bibitem{ka} T. Kato, {\it Quasilinear equations of evolution with applications to
partial differential equations}, Springer Lecture Notes {\bf 448},
1975, 25-70.

\bibitem{k} S. Klainerman, {\it Global existence for nonlinear wave
equations}, Comm. Pure Appl. Math. {\bf 33} (1980), 43-101.

\bibitem{kl}
De-Xing Kong and Kefeng Liu, {\it Wave character of metrics and
hyperbolic geometric flow},
http://www.cms.zju.edu.cn/UploadFiles/AttachFiles/200682885946597.pdf.

\bibitem{p} R. Penrose, {\it Gravitational collapse and space-time singularities},
Phys. Rev. Lett. {\bf 14} (1965), 57-59.

\bibitem{s1} J. Shatah and M. Struwe, Geometric Wave Equations, Courant
Lecture Notes in Mathematics {\bf 2}, New York University, Courant
Institute of Mathematical Sciences, New York; American Mathematical
Society, Providence, RI, 1998.

\bibitem{s2}  J. Shatah and M. Struwe, {\it Regularity results for nonlinear wave equations},
Ann. of Math. {\bf 138} (1993), 503-518.

\bibitem{s3} J. Shatah and A. Tahvildar-Zadeh, {\it Regularity of harmonic maps from the
Minkowski space into rotationally symmetric manifolds}, Comm. Pure
Appl. Math. {\bf 45} (1992), 947-971.

\bibitem{shu} Fu-Wen Shu and You-Gen Shen, {\it Geometric flows and black holes}, arXiv: gr-qc/0610030.

\end{thebibliography}
\end{document}